\documentclass{article}
   \usepackage{amssymb} 
\newcommand{\g}{\mbox{$\bf g$}}

\newcommand{\h}{\mbox{\textbf{h}}}
\newcommand{\n}{\mbox{\textbf{n}}}
\newcommand{\np}{\mbox{$\textbf{n}^+$}}
\newcommand{\nm}{\mbox{$\textbf{n}^-$}}
\newcommand{\al}{\alpha}
\newcommand{\eps}{\epsilon}
\newcommand{\la}{\lambda}
\newcommand{\La}{\Lambda}
\newcommand{\de}{\delta}
\newcommand{\Th}{\Theta}

\newcommand{\mb}{\mbox}

\newcommand{\Mklz}[2]{\left\{\left.\;#1\;\right|\; #2\;\right\}}
\newcommand{\W}{\mbox{$\Delta$}}
\newcommand{\nW}{\mbox{$\Delta^-$}}
\newcommand{\pW}{\mbox{$\Delta^+$}}
\newcommand{\rW}{\mbox{$\Delta_{re}$}}
\newcommand{\nrW}{\mbox{$\Delta_{re}^-$}}
\newcommand{\prW}{\mbox{$\Delta_{re}^+$}}
\newcommand{\iW}{\mbox{$\Delta_{im}$}}

\newcommand{\F}{\mathbb{F}}

\newcommand{\C}{\mathbb{C}}
\newcommand{\N}{\mathbb{N}}
\newcommand{\Nn}{\mathbb{N}_0}
\newcommand{\Q}{\mathbb{Q}}
\newcommand{\Qp}{\mathbb{Q}^+}

\newcommand{\R}{\mathbb{R}}

\newcommand{\Z}{\mathbb{Z}}
\newcommand{\We}{\mbox{$\cal W$}}

\newcommand{\Cq}{\mbox{$\overline{C}$}}

\newcommand{\RkX}{\mbox{${\cal R}(X)$}}
\newcommand{\iB}[2]{\left(#1\mid#2\right)}

\newcommand{\kB}[2]{\left\langle \left\langle #1\mid #2
                       \right\rangle \right\rangle}
                       
\newcommand{\kBl}{\left\langle \left\langle \;\ \mid \;\ \right\rangle 
                       \right\rangle}
\newcommand{\TD}{\mbox{$\widehat{T}$}}
\newcommand{\ND}{\mbox{$\widehat{N}$}}
\newcommand{\GD}{\mbox{$\widehat{G}$}}
\newcommand{\GfD}{\mbox{$\widehat{G_f}$}}
\newcommand{\WeD}{\mbox{$\widehat{\We}$}}
\newcommand{\ve}[1]{\mbox{$\varepsilon\left(#1\right)$}} 
\newcommand{\FK}[1]{\mbox{$\F\,[#1]$}}
\newcommand{\CK}[1]{\mbox{${\mathbb C}\,[#1]$}}

\newcommand{\Proof}{\mbox{\bf Proof: }}
\newcommand{\End}{\mb{}\hfill\mb{$\square$}\\}

\newcommand{\Spm}{\mbox{Specm\,}}

\newcommand{\ti}{\tilde}
\newcommand{\res}[1]{\!\mid_{#1}}

\newcommand{\tr}{{\,\bf \diamond\,}}
\begin{document}
\newtheorem{Theorem}{Theorem}[section]
\newtheorem{Def}[Theorem]{Definition}
\newtheorem{Prop}[Theorem]{Proposition}
\newtheorem{Prop+Def}[Theorem]{Proposition+Definition}
\newtheorem{Cor}[Theorem]{Corollary}
\newtheorem{Rem}[Theorem]{Remark}
\newtheorem{Rems}[Theorem]{Remarks}
\newtheorem{Lemma}[Theorem]{Lemma}
%
%
\title{Extending the Bruhat order and the length function from the Weyl group to the Weyl monoid}
\author{Claus Mokler\\\\ Universit\"at Wuppertal, Fachbereich Mathematik\\  Gau\ss stra\ss e 20\\ D-42097 Wuppertal, Germany\vspace*{1ex}\\ 
          mokler@math.uni-wuppertal.de}
\date{}
\maketitle
\begin{abstract}\noindent
For a symmetrizable Kac-Moody algebra the category of admissible representations is an analogue of the 
category of finite dimensional representations of a semisimple Lie algebra. 
The monoid associated to this category and the category of restricted duals by a generalized Tannaka-Krein 
reconstruction contains the Kac-Moody group as open dense unit group and has similar properties as a 
reductive algebraic monoid. In particular there are Bruhat and Birkhoff decompositions, the Weyl group 
replaced by the Weyl monoid, \cite{M1}.\\
We determine the closure relations of the Bruhat and Birkhoff cells, which give extensions of the Bruhat order 
from the Weyl group to the Weyl monoid. We show that the Bruhat and Birkhoff cells are irreducible and 
principal open in their closures. We give product decompositions of the Bruhat and Birkhoff cells.
We define extended length functions, which are compatible with the extended Bruhat orders. We show a 
generalization of some of the Tits axioms for twinned BN-pairs.
\end{abstract}
\vspace*{1ex}
\section*{Introduction}
The Kac-Moody group $G$ constructed in \cite{KP1} by V. Kac and D. Peterson is a group analogue 
of a semisimple, simply connected algebraic group. In particular there are Bruhat and Birkhoff decompositions
\begin{eqnarray*}
   G\;\,=\;\,\dot{\bigcup_{w\in{\cal W}}} B^\eps w B^\delta && 
                               (\eps,\delta)\;\in\; \{\,(+,+),\,(-,-),\,(-,+),\,(+,-)\,\}\;\;,
\end{eqnarray*}  
and also certain multiplicative decompositions of the Bruhat and Birkhoff cells like $BwB= U_w\cdot wT\cdot U$ and 
$B^-wB= (U^w)^- \cdot wT\cdot U$. (Here $T$ denotes a maximal torus, $\We$ the Weyl group, and $B=B^+$, $B^-$ denote 
opposite Borel subgroups containing $T$.)
V. Kac and D. Peterson equipped in \cite{KP2} a symmetrizable Kac-Moody group with a coordinate ring, the algebra 
of strongly regular functions $\CK{G}$. This coordinate ring has many properties in common with the coordinate ring 
of a semisimple, simply connected algebraic group. It is an integrally closed domain, even a unique 
fac\-tor\-iza\-tion domain. It admits a Peter and Weyl theorem, i.e.,
\begin{eqnarray*}
  \CK{G} \;\cong \; \bigoplus_{\La\in P^+} L^*(\La)\otimes L(\La)
\end{eqnarray*}
as $G\times G$-modules. (Here $L(\La)$ denotes an irreducible highest weight module with highest weight $\La$, 
$L^*(\La)$ its restricted dual, and $P^+$ the set of dominant weights.) As a difference, in the non-classical case, 
the multiplication map and the inverse map of $G$ do not induce comorphisms. V. Kac and D. Peterson showed that the 
Zariski closures of the following Bruhat and Birkhoff cells are obtained similarly to the classical case, using 
the Bruhat order of the Weyl group $\We$:
\begin{eqnarray*}
  \overline{B^\eps w B^\eps}\; =\; \bigcup_{w'\leq w }B^\eps w' B^\eps \quad\mb{where}\quad\eps\in\{+,-\}, 
   \quad\quad \overline{B^- w B}\; =\; \bigcup_{w'\geq w }B^-w' B\;\;.
\end{eqnarray*}
In \cite{M1} we determined by a generalized Tannaka-Krein reconstruction the monoid with coordinate ring 
($\GD$, $\CK{\GD}$) associated to a natural category determined by the modules $L(\La)$, $\La\in P^+$, and a 
category of duals determined by $L^*(\La)$, $\La\in P^+$. For its history in connection with V. Kac, D. Peterson, 
and P. Slodowy please look at the introduction of \cite{M1}.
The monoid $\GD$ contains the Kac-Moody group $G$ as open dense unit group. In particular its coordinate ring 
$\CK{\GD}$ is isomorphic to the algebra of strongly regular functions $\CK{G}$ by the restriction map.\\
This monoid is a purely infinite-dimensional phenomenon. In the classical case it reduces to a semisimple simply 
connected algebraic group. It is a proper analogue of such a group. For generalizing some results of 
classical invariant theory it should be more fundamental than the Kac-Moody group itself.\\
In \cite{M1} we also showed that the monoid $\GD$ has similar structural properties as a normal reductive algebraic 
monoid. In particular there are Bruhat and Birkhoff decompositions
\begin{eqnarray*}
   \GD\;\,=\;\,\dot{\bigcup_{\hat{w}\in \widehat{\cal W}}}B^\eps \hat{w} B^\delta && 
                               (\eps,\delta)\;\in\; \{\,(+,+),\,(-,-),\,(-,+),\,(+,-)\,\}\;\;,
\end{eqnarray*}  
the Weyl group replaced by the Weyl monoid $\WeD$, which is an analogue of a Renner monoid. It contains the 
Weyl group as unit group and its idempotents correspond bijectively to the faces of the Tits cone. \\ 
In \cite{M2} we determined and investigated the $\C$-valued points of $\CK{\GD}$. Identifying the elements of 
$\GD$ with their evaluation morphisms, $\GD$ embeds in the set of $\C$-valued points of $\CK{\GD}$.  
The Bruhat decompositions of $\GD$ do not extend, but one of the Birkhoff decompositions of $\GD$ extends to a 
decomposition of the $\C$-valued points.\\ 
In \cite{M1}, \cite{M2} the Bruhat and Birkhoff cells have not been investigated further. In particular their 
closure relations, which determine extensions of the Bruhat order from the Weyl group to the Weyl monoid, have 
not been determined.\vspace*{1ex}\\
For reductive algebraic monoids these questions have already been investigated:\\
L. E. Renner studied in \cite{Re1}, \cite{Re2} the closure relation of the Bruhat cells of a reductive algebraic 
monoid. Transferred to the Renner monoid he called this order the Bruhat-Chevalley order. He showed that all maximal 
chains between two elements have the same length. He also showed a monoid version of one of the Tits axioms for 
$BN$-pairs. An algebraic description of the Bruhat-Chevalley order has been obtained by E. A. Pennel, M. S. Putcha, 
and R. E. Renner in \cite{PePuRe}. An investigation of the lexicographic shellability has been started by 
M. S. Putcha in \cite{Pu}.\\  
A length function for the Renner monoid of matrices over a finite field has first been introduced in \cite{So} 
by L. Solomon. He also showed that this length function fits to the monoid 
generalization of one of the Tits axiom mentioned above.
A length function for the Renner monoid of an arbitrary reductive algebraic monoid has been introduced and 
investigated by L. E. Renner in \cite{Re2} by a different approach. It has the property that the 
length of a maximal chain of the Bruhat-Chevalley order, which is contained in an orbit of the action of the 
product of the Weyl group on the Renner monoid, is given by the difference of the length of the maximal and 
the minimal element of the chain. He also showed that this length function fits to his monoid generalization of one 
of the Tits axioms. A subadditivity property has been shown by E. A. Pennel, M. S. Putcha, and R. E. Renner in 
\cite{PePuRe}.\\
L. E. Renner gave in \cite{Re3} a product decomposition of a Bruhat cell of the wonderful compactification of a 
semisimple algebraic group.\\
Note also that many of these results just mentioned induce, by using the longest element of the Weyl group, 
corresponding results for the Birkhoff cells.\vspace*{1ex}\\
In this article we show that some of these results are also valid for the monoid $\GD$ and the spectrum 
$\Spm\CK{\GD}$:
In Section \ref{SBO} we obtain a similar algebraic description of the closure relations of the Bruhat and Birkhoff 
cells of $\GD$, and of the Birkhoff cells of $\Spm\CK{\GD}$. We show that these relations are actually order 
relations, which is not automatic in our infinite dimensional situation. Transferred to the Weyl monoid we call 
these orders the extended Bruhat and Birkhoff orders.\\
We show that the Bruhat and Birkhoff cells are irreducible and principal open in their closures. We equip the 
Bruhat and Birkhoff cells with their coordinate rings as principal open sets and
give product decompositions of these cells.\\
In Section \ref{SL} we extend the length function on the Weyl group to functions on the Weyl 
monoid, compatible with the extended Bruhat orders. We show a generalization of some of the Tits axiom for groups with 
twinned $BN$-pairs.\vspace*{1ex}\\
To obtain these results is often more arduous as for reductive algebraic monoids, and we have to use different 
methods for the following reasons:
Most of the theorems of algebraic geometry, which are used to investigate algebraic groups and monoids, break 
down for the infinite-dimensional varieties we have to use. In the non-classical case the multiplication map of 
$\GD$ is not a morphism, only the left and right multiplications with elements of $\GD$ are morphisms. The Weyl 
group and the Weyl monoid have infinitely many elements. In particular there is no longest element of the Weyl 
group, and there are infinitely many Bruhat and Birkhoff cells.\vspace*{1ex}\\
For our algebraic geometric investigations we mainly use the following two aids:
The explicit description of the coordinate ring of $\GD$ by the Peter and Weyl theorem. This allows to use 
properties of the action of parts of $\GD$ on the irreducible highest weight representations $L(\La)$, 
$\La\in P^+$, to obtain algebraic geometric results. As an advantage our morphisms, closures, principal open 
sets are described very explicitely by using matrix coefficients.
The monoid $\GD$ does not act on the flag varieties, but it acts on the corresponding affine cones by morphisms. 
This allows to make use of the combinatorial properties of these cones described by V. Kac and D. Peterson in 
\cite{KP1}.
%
%
%
%
%
%
%
%
%
%
%
%
%
%
%
%
%
\tableofcontents\mb{}\\
%
%
%
%
%
%
%
%
\section{Preliminaries}
%
%
%
In this section we collect some basic facts about Kac-Moody algebras, minimal and formal Kac-Moody groups, and the 
cor\-res\-pond\-ing monoid completions, which are used later.
One aim is to introduce our notation. Another aim is to put these things, which are scattered between many articles and 
books, on equal footing appropriate for our goals.\\
All the material stated in this section about Kac-Moody algebras can be found in the books \cite{K} (most results also valid 
for a field of characteristic zero with the same proofs), \cite{MoPi}. The facts about the minimal Kac-Moody group can be found in 
\cite{KP1}, \cite{KP3}, \cite{MoPi}, about the formal Kac-Moody group in \cite{Sl1}, about the Kostant cones in \cite{KP1}, \cite{MoPi}, 
about the algebra of strongly regular functions in \cite{KP2}, about the faces of the Tits cone in \cite{Loo}, \cite{Sl1}, \cite{M1}.
The facts about the monoid completion of the minimal Kac-Moody group can be found in \cite{M1}, about the spectrum of $\F$-valued points of its 
coordinate ring in \cite{M2}.\vspace*{1ex}\\ 
We denote by $\N=\Z^+$, $\Q^+$, resp. $\R^+$ the sets of strictly positive numbers of $\Z$, $\Q$, resp. $\R\,$,
and the sets $\N_0=\Z^+_0$, $\Q^+_0$, $\R^+_0$ contain, in addition, the zero.
In the whole paper, $\F$ is a field of characteristic 0 and $\F^\times$ its group of units.\vspace*{1ex}\\
%
%
%
%
%
{\bf Generalized Cartan matrices:}
%
%
%
Starting point for the construction of a Kac-Moody algebra and its associated simply connected minimal and formal Kac-Moody groups is a 
{\it generalized Cartan matrix}, which is a matrix $A=(a_{ij})\in M_{n}(\Z)$ with $a_{ii}=2$, $a_{ij}\leq 0$ for all $i\neq j$, and $a_{ij}=0$ 
if and only if $a_{ji}=0$. Denote by $l$ the rank of $A$, and set $I:=\{1,2,\ldots, n\}$.\\ 
For the properties of the generalized Cartan matrices, in particular their classification, we refer to the book \cite{K}. In this paper 
we assume $A$ to be symmetrizable.\vspace*{1ex}\\
%
%
%
{\bf Realizations:}
%
%
%
A {\it simply connected minimal free realization} of $A$ consists of dual free  $\mathbb{Z}$-modules $H$, $P$ of rank $2n-l$, 
and linear independent sets $\Pi^\vee =\{h_1,\ldots, h_n\}\subseteq H $, $\Pi=\{\al_1,\ldots,\al_n\}\subseteq P$ such that 
$\al_i(h_j)=a_{ji}\,$, $i,j=1,\dots, n$. Furthermore there exist (non-uniquely determined) fundamental do\-mi\-nant weights
$\La_1,\ldots, \La_n\in P$ such that $\La_i(h_j) =\de_{ij}$, $i,j=1,\ldots, n$.
$P$ is called the {\it weight lattice}, and $Q:=\Z\mb{-span}\Mklz{\al_i\,}{\,i\in I}$ the {\it root lattice}.
Set $Q^\pm_0:=\Z^\pm_0\mb{-span}\Mklz{\al_i\,}{\,i\in I}$, and $Q^\pm:=Q^\pm_0\setminus\{0\}$.\\
We fix a system of fundamental dominant weights $\La_1,\ldots, \La_n$, and extend $h_1,\ldots, h_n\in H$, 
$\La_1,\ldots, \La_n\in P$ to a pair of dual bases $h_1,\ldots, h_{2n-l}\in H$, $\La_1,\ldots, \La_{2n-l}\in P$. We
set $H_{rest}:=\Z \mb{-}span\Mklz{h_i}{i=n+1,\ldots,2n-l}$.\vspace*{1ex}\\
%
%
%
%
%
%
%
{\bf The Weyl group:}
%
%
%
%
Identify $H$ and $P$ with the cor\-res\-pond\-ing sublattices of the following vector spaces over $\F\,$:
\begin{eqnarray*}
   \h  \;\,:=\;\,  \h_\F    \;\,:=\;\,   H \otimes_{\mathbb Z} \F   &\quad,\quad &
   \h^* \;\,:=\;\, \h^*_\F  \;\,:=\;\,  P \otimes_{\mathbb Z} \F\;\;.
\end{eqnarray*}
$\h^*$ is interpreted as the dual of $\h$. Order the elements of $\h^*$ by $\la\leq\la'$ if and only if $\la'-\la\in Q_0^+$.\\ 
Choose a symmetric matrix $B\in M_n (\Q)$ and a diagonal matrix 
$D=\mb{diag}(\eps_1,\ldots, \eps_n)\,$, $\eps_1,\ldots,\eps_n\in\Qp$, such that $A=DB$.
Define a nondegenerate symmetric bilinear form on $\h$ by:
\begin{eqnarray*}
 \iB{h_i}{h} \;=\; \iB{h}{h_i} \;:=\; \al_i(h)\,\eps_i  &\qquad & i\in I\,,
 \quad h\in\h\;,\\
 \iB{h'}{h''} \;:=\; 0 \qquad\qquad &\qquad & h',\,h'' \in \h_{rest}:=H_{rest}\otimes\F \;.
\end{eqnarray*}
Denote the induced nondegenerate symmetric form on $\h^*$ also by $\iB{\;}{\;}$.\vspace*{1ex}\\
The {\it Weyl group} $\We=\We(A)$ is the Coxeter group with generators $\sigma_i\,$, $i\in I$, and
relations 
\begin{eqnarray*}
       \sigma_i^2 \;=\; 1 \qquad (i\in I)    \;\:&\;,\;&\;\:
       {(\sigma_i\sigma_j)}^{m_{ij}} \;=\; 1 \qquad (i,j\in I,\,i\ne j)\;\;.
\end{eqnarray*}
The $m_{ij}$ are given by:
    $\quad \begin{tabular}{c|ccccc}
      $a_{ij}a_{ji}$ & 0 & 1 & 2 & 3 &  $\geq$ 4  \\[0.5ex] \hline 
         $m_{ij}$    & 2 & 3 & 4 & 6 &  no relation between $\sigma_i$ and $\sigma_j$ 
     \end{tabular}$\vspace*{1ex}\\
The Weyl group $\We$ acts faithfully and contragrediently on $\h$ and $\h^*$ by 
\begin{eqnarray*}
  \sigma_i h \;:=\; h - \al_i\,(h) h_i \quad,\quad
  \sigma_i \la \;:=\; \la - \la(h_i)\,\al_i  \quad,\quad i\in I,\quad h\in \h, \quad \la\in \h^*\;\;,
\end{eqnarray*}
leaving the lattices $H$, $Q$, $P$, and the forms invariant.
$\Delta_{re}:=\We\Mklz{\al_i}{i\in I}$ is called the set of {\it real roots}, and 
$\Delta_{re}^\vee:=\We\Mklz{h_i}{i\in I}$ the set of {\it real 
coroots}. The map $\al_i\mapsto h_i\,$, $i\in I$, can be extended to a $\We$-equivariant bijection
$\al \mapsto  h_\al$.\vspace*{1ex}\\
%
%
%
{\bf The Tits cone and its faces:}
To illustrate the action of $\We$ on $\h^*_\R$ geometrically, for $J\subseteq I$ define 
\begin{eqnarray*}
  F_J &:=& \Mklz{\la\in\h^*_\R}{\la(h_i)\,=\,0 \;\mb{ for }\; i\in J\,,\;\;\;
  \la(h_i)\,>\,0\; \mb{ for }\;i\in I\setminus J}\;\;,\\
  \overline{F_J} &:=& \Mklz{\la\in\h^*_\R}{\la(h_i)\,=\,0\; \mb{ for }\;i\in J
  \,,\;\;\;\la(h_i)\,\geq\,0\; \mb{ for }\;i\in I\setminus J}\;\;. 
\end{eqnarray*}
$\overline{F_J}$ is a finitely generated convex cone with relative interior $F_J$. The parabolic subgroup $\We_J$ of $\We$ is the stabilizer of every 
element $\la\in F_J$. For $\sigma\in\We$ call $\sigma F_J$ a {\it facet} of {\it type} $J$.\\   
The {\it fundamental chamber} $\overline{C} \,:=\, \Mklz{\la\in\h^*_\R}{\la(h_i)\,\geq\,0\; \mb{ for }\;i\in I}$ is a fundamental 
region for the action of $\We$ on the convex cone $X:=\We\,\overline{C}\,$, which is called the {\it Tits cone}. The partition
$\overline{C}=\dot{\bigcup}_{J\subseteq I} F_J$ induces a $\We$-invariant partition of $X$ into facets.\vspace*{1ex}\\ 
We denote the set of faces of the Tits cone $X$ by $\RkX$. These faces can be described as follows: A set $\Th\subseteq I$ is called {\it special}, 
if either $\Th=\emptyset$, or else all connected components of the generalized Cartan submatrix $(a_{ij})_{i,j\in\Th}$ are of non-finite type. 
Set $\Th^\bot:=\Mklz{i\in I}{ a_{ij}=0 \mb{ for all } j\in\Th }$. Every face of the Tits cone $X$ is $\We$-conjugate to exactly one of the faces
\begin{eqnarray*}
  R(\Th) \;\,:=\;\, X\,\cap\,\Mklz{\la\in \h_{\R}^*}{\la(h_i)=0 \mb{ for all } i\in\Th} 
          \;\,=\;\,\We_{\Th^\bot} \overline{F}_\Th&\;,\; & \Th\;\mb{ special}\;\,.
\end{eqnarray*}
The parabolic subgroup $\We_\Th$ is the pointwise stabilizer of $R(\Th)$, and the parabolic subgroup 
$\We_{\Th\cup\Th^\bot}$ is the stabilizer of the set $R(\Th)$ as a whole. \\
The relative interior of $R(\Th)$ is given by the union of the facets $\sigma F_{\Th\cup\Th_f}$, where $\sigma\in\We_{\Th^\bot}$, and 
$\Th_f$ is a subset of $\Th^\bot$, which is either empty, or else for which all connected components of $(a_{ij})_{i,j\in\Th_f}$ are of 
finite type.\vspace*{1ex}\\
%
%
%
{\bf The Weyl monoid:}
The Weyl group acts on the monoid (\,$\RkX\,,\,\cap$\,). The semidirect product $\RkX\rtimes\We$ consists of the set $\RkX\times\We$ with the 
structure of a monoid given by 
\begin{eqnarray*}
  (R,\sigma)\cdot(S,\tau) &:=& (R\cap\sigma S,\sigma \tau)\;\;.
\end{eqnarray*}
For $R\in\RkX$ let $Z_{\cal W}(R):=\Mklz{\sigma\in \We}{\sigma\la=\la \mb{ for all }\la\in R}$ be the pointwise stabilizer of $R$. The Weyl monoid 
$\WeD$ is defined as the monoid $\RkX\rtimes\We$ factored by the congruence relation
\begin{eqnarray*}
   (R,\sigma) \sim (R',\sigma') &:\iff & R\:=\:R'\quad\mb{and}\quad\sigma'\sigma^{-1}\in Z_{\cal W}(R)\;\;.
\end{eqnarray*}
We denote the congruence class of $(R,\sigma)$ by $\ve{R}\sigma$.\\
Assigning to $\sigma\in\We$ the element $\sigma:=\ve{X}\sigma \in\WeD$, the Weyl group $\We$ identifies with the 
unit group of $\WeD$. It acts in the obvious way on $\WeD$, the partition of $\WeD$ into $\We\times\We$-orbits 
given by
\begin{eqnarray}\label{WWorbits}
  \WeD &=& \dot{\bigcup_{\Th \;special}} \,\We\,\ve{R(\Th)}\,\We\;\;.
\end{eqnarray} 
Assigning to $R\in \RkX$ the element $\ve{R}:=\ve{R}1\in\WeD$, the monoid $(\,\RkX\,,\,\cap \,)$ embeds into $\WeD$. 
Its image are the idempotents of $\WeD$. 
By this map, the action of $\We$ on $\RkX$ identifies with the restricted conjugation action of $\We$ on $\WeD$, 
i.e., 
\begin{eqnarray}\label{wRwinv}
  w e(R)w^{-1}=e(wR)\;,\quad\quad  w\in\We\;,\quad R\in\RkX\;.
\end{eqnarray}
Recall that a monoid $M$ is called an {\it inverse monoid}, if for every element $m\in M$ there exists a unique 
element $m^{inv}\in M$, such that $m m^{inv} m = m$ and $m^{inv}m m^{inv}= m^{inv}$. The {\it inverse map} 
$\mb{}^{inv}:M\to M$ is an involution extending the inverse map of the unit group.
The Weyl monoid $\WeD$ is an inverse monoid with inverse map $\mb{}^{inv}:\WeD\to\WeD$ given by 
$(\sigma e(R))^{inv}=e(R)\sigma^{-1}$, $\sigma\in \We$, and $R$ a face of $X$.\\
The following formulas, which follow immediately from the definition of $\WeD$, will be used later at many 
places: Let $w\in\We$ and $\Th$ be special. Then
\begin{eqnarray}
   w\ve{R(\Th)}w^{-1}\;=\;\ve{R(\Th)} \;\,\iff\;\, w\in\We_{\Th\cup\Th^\bot} \;\;,\label{WeDF1}\hspace{5em}\\
   w\ve{R(\Th)}\;=\;\ve{R(\Th)}  \;\,\iff\;\,      \ve{R(\Th)}w\;=\;\ve{R(\Th)} \;\,\iff\;\, w\in\We_\Th\;\;.
   \label{WeDF2}
\end{eqnarray}
%
%
%
%
{\bf The Kac-Moody algebra:} 
%
%
%
%
The {\it Kac-Moody algebra} $\g=\g(A)$ is the Lie algebra over $\F$ generated by the abelian Lie algebra $\h$ and
$2n$ elements $e_i,f_i$, ($i\in I$), with the following relations, which 
hold for any $i,j \in I$, $h \in \h\,$: 
  \begin{eqnarray*}
    \left[ e_i,f_j \right] \,=\,   \delta_{ij} h_i   \;\;,\;\;   
    \left[ h,e_i \right]   \,=\,   \al_i(h) e_i    \;\;,\;\;  
    \left[ h,f_i \right]   \,=\,  -\al_i(h) f_i \;\;,    \\
    \left(ad\,e_i\right)^{1-a_{ij}}e_j  \,=\, \left(ad\,f_i\right)^{1-a_{ij}}f_j  \,=0  \qquad (i\neq j)\;\;.
  \end{eqnarray*}
The {\it Chevalley involution} $*:\g\to\g$ is the involutive anti-automorphism determined by $ e_i^*=f_i$, $ f_i^*=e_i$, $h^*=h$, 
($i\in I$, $h\in \h$).\\     
The nondegenerate symmetric bilinear form \mb{( $|$ )} on $\h$ extends uniquely to a non\-de\-ge\-ner\-ate symmetric invariant 
bilinear form \mb{( $|$ )} on $\g$.\\
We have the {\it root space decomposition}
\begin{eqnarray*}
  \g=\bigoplus_{\al \in {\bf h}^*}\g_{\al} \quad \mb{where} \quad 
  \g_\al := \Mklz{x\in \g}{[h,x]=\al(h)\,x\;\mb{ for all }\;h\in \h}\;\;.           
\end{eqnarray*}
In particular $\g_0 = \h$, $\g_{\al_i}=\F e_i$, and $\g_{-\al_i}=\F f_i$, $i\in I $.\\
The set of roots $\W:=\Mklz{\al\in\h^*\setminus\{0\}}{\g_\al\ne \{0\}}$ is 
invariant under the Weyl group, $\W=-\W$, and $\W$ spans the root lattice $Q\,$. We have $\rW\subseteq \W\,$, and 
$\iW:=\W\setminus\rW$ is called the set of {\it imaginary roots}.\\
$\W$, $\rW$, and $\iW$ decompose into the disjoint union of the sets of 
{\it positive} and {\it negative} roots $\W^\pm:=\W\cap Q^\pm$, $\rW^\pm:=\rW\cap Q^\pm$, $\iW^\pm:=\iW\cap Q^\pm$.\\ 
There is the {\it triangular decomposition} $\g = \nm \oplus \h \oplus \np $, where $\n^\pm:=\bigoplus_{\al\in {\Delta}^\pm} \g_\al $.\vspace*{1ex}\\        
%
%
%
%
%
{\bf Irreducible highest weight representations:}
%
%
%
%
For every $\La\in\h^*$ there exists, unique up to isomorphism, an irreducible representation $(L(\La),\pi_\La)$
of $\g$ with highest weight $\La$. It is $\h$-dia\-gon\-al\-iz\-able, and we denote its set 
of weights by $P(\La)$.
Any such representation carries a nondegenerate symmetric bilinear form $\kBl:L(\La)\times L(\La)\to\F$
which is contravariant, i.e., $\kB{v}{xw}=\kB{x^*v}{w}$ for all $v,w\in L(\La)$, $x\in\g$. This form is unique 
up to a nonzero multiplicative scalar.\vspace*{1ex}\\  
%
%
%
%
{\bf The category ${\cal O}_{adm}$:}
%
%
%
The category $\cal O$ is defined as follows: Its objects are the $\g$-modules $V$, which have the properties:\\
(1) $V$ is $\h$-diagonalizable with finite dimensional weight spaces.\\
(2) There exist finitely many elements $\la_1,\,\ldots,\,\la_m\in\h^*$, such that the set of weights $P(V)$ of $V$ is contained in the union 
$\bigcup_{1=1}^m \Mklz{\la\in\h^*}{\la\leq \la_i}$.\\
The morphisms of $\cal O$ are the morphisms of $\g$-modules.\vspace*{1ex}\\ 
Call a $\g$-module $V$ {\it admissible}, if $V$ is $\h$-diagonalizable with set of weights $P(V)\subseteq P$ and 
the elements of $\g_\al$ act locally nilpotent on $V$ for all $\al\in\rW$.
(If the generalized Cartan matrix is degenerate, then admissible is slightly stronger than {\it integrable}, 
which means $V$ is $\h$-diagonalizable and the elements of $\g_\al$ act locally nilpotent on $V$ for all 
$\al\in\rW$.)
Examples of admissible representations are the adjoint representation ($\g\,$, $ad\,$), and the irreducible highest 
weight representations ($L(\La)$, $\pi_\La$), $\La\in P^+:=P\cap\Cq$.\vspace*{1ex}\\
We denote by ${\cal O}_{adm}$ the full subcategory of the category $\cal O$, whose objects are admissible modules. This category generalizes the 
category of finite dimensional representations of a semisimple Lie algebra, keeping the complete reducibility theorem. Every object of 
${\cal O}_{adm}$ is isomorphic to a direct sum of the admissible irreducible highest weight modules $L(\La)$, $\La\in P^+$.
The set of weights of a module of ${\cal O}_{adm}$ is contained in $X\cap P$ because of 
$\bigcup_{\La\in P^+} P(\La) = X\cap P$.\vspace*{1ex}\\
%
%
%
%
%
%
{\bf The minimal and formal Kac-Moody groups $G$ and $G_f$, the monoids $\GD$ and $\GfD$:}
%
%
The monoid $\GD$, can be obtained by a Tannaka-Krein reconstruction from the category ${\cal O}_{adm}$ and its corresponding category of restricted 
duals, compare \cite{M1}, Section 4. It can also be characterized as follows:\\
{\bf (a)} The monoid $\GD$ acts on every module of ${\cal O}_{adm}$. Two elements $\hat{g},\hat{g}'\in \GD$ are equal if and only if  
for all modules $V$ of ${\cal O}_{adm}$ and for all $v\in V$ we have $\hat{g}v=\hat{g}'v$.\\ 
{\bf (b)} There are the following elements of $\GD$ acting on the modules of ${\cal O}_{adm}$ in a particular way:\\
(1) For every $h\in H$, $s\in\F^\times$ there exists an element $t_h(s)\in \GD$, such that for every module $V$ of ${\cal O}_{adm}$
we have
\begin{eqnarray*}
  t_h(s)v_\la &=& s^{\la(h)}v_\la \quad,\quad v_\la\in V_\la\;\,,\;\,\la\in P(V)\;.
\end{eqnarray*}
(2) For every $x\in\g_\al$, $\al\in\rW$, there exists an element $\exp(x)\in \GD$, such that for 
every module $V$ of ${\cal O}_{adm}$ we have
\begin{eqnarray*}
  \exp(x)v &=& \exp(\pi(x)) v \quad,\quad v\in V\;.
\end{eqnarray*}
(3) For every face $R$ of the Tits cone there exists an element $e(R)\in \GD$, such that for every module $V$ of ${\cal O}_{adm}$ we have
\begin{eqnarray*}
  e(R)v_\la &=& \left\{ \begin{array}{ccc}
    v_\la &\;&\la\in R\\
      0   &\;& \la\in X\setminus R
  \end{array}\right. \quad,\quad \;v_\la\in V_\la\,,\; \;\la\in P(V)\;\;.
\end{eqnarray*}
$\GD$ is generated by the elements of (1), (2), and (3).\vspace*{1ex}\\
The unit group $G$ of $\GD$ is generated by the elements of (1) and (2). It is the {\it minimal Kac-Moody group}, which we call {\it Kac-Moody group} 
for short.\vspace*{1ex}\\
The {\it Chevalley involution} $*:\,\GD \to \GD$ is the involutive anti-isomorphism determined by $\exp(x_\al)^*:=\exp(x_\al^*)$, 
$t^*:= t$, $e(R)^* := e(R)$, where $x_\al\in\g_\al$, $\al\in\rW$, $t\in T$, and $R\in\RkX$.
It is compatible with any nondegenerate symmetric contravariant form $\kBl$ on any module $V$ of ${\cal O}_{adm}$, i.e., $\kB{xv}{w}=\kB{v}{x^*w}$, 
$v,w\in V$, $x\in\GD$.\vspace*{1ex}\\
In this paper we are interested in the monoid $\GD$ and the spectrum of $\F$-valued points of its coordinate ring, which will be defined soon. To 
describe the $\F$-valued points we need a second monoid $\GfD$, extending $\GD$, which we define already now:
Set $\n_f:=\prod_{\al\in \Delta^+}\g_\al$ and $\g_f:=\n^-\oplus\h\oplus\n_f$. The Lie bracket of $\g$ extends in the obvious
way to a Lie bracket of $\g_f$. Every $\g$-module of ${\cal O}_{adm}$ can be extended to a $\g_f$-module.
The monoid $\GfD$ can be characterized as follows:\\
{\bf (a)} The monoid $\GfD$ acts on every module of ${\cal O}_{adm}$. Two elements $\hat{g},\hat{g}'\in \GD$ are equal if and only if  
for all modules $V$ of ${\cal O}_{adm}$, and for all $v\in V$, we have $\hat{g}v=\hat{g}'v$.\\ 
{\bf (b)} $\GfD$ extends $\GD$ and it contains the following elements:\\
(4) For every $x\in\n_f$ there exists an element $\exp(x)\in \GfD$, such that for any admissible representation $(V,\pi)$ we have
\begin{eqnarray*}
  \exp(x)v &=& \exp(\pi(x)) v \quad,\quad v\in V\;\;.
\end{eqnarray*}
$\GfD$ is generated by $\GD$ and the elements of (4).\vspace*{1ex}\\
The unit group $G_f$ of $\GfD$ is generated by $G$ and the elements of (4). It is the {\it formal Kac-Moody group}.\vspace*{1ex}\\
Note that the groups $G$, $G_f$ as well as the monoid $\GD$, $\GfD$ act faithfully on the sum $\bigoplus_{\La\in P^+}L(\La)$.\vspace*{1ex}\\  
The groups $G$ and $G_f$ have the following important structural properties:\vspace*{0.5ex}\\
{\bf (a)} The elements of (1) induce an embedding of the torus $H\otimes_\Z\F^\times$ into $G\subseteq G_f$. Its image is denoted by $T$.
For $\al\in\rW$ the elements of (2) induce an embedding of $(\g_\al,+)$ into $G\subseteq G_f$. Its image $U_\al$ is 
called the {\it root group} belonging to $\al$.\\
Let $\al\in \prW$ and $x_\al\in\g_{\al}$, $x_{-\al}\in\g_{-\al}$ such that $[x_\al,x_{-\al}]=h_\al$. There exists an 
injective homomorphism of groups $\phi_\al:\,\mb{SL}(2,\F) \to G$ with 
\begin{eqnarray*}
 \phi_\al\left(\begin{array}{cc}
  1 & s\\
  0 & 1
  \end{array}\right) \;:=\; \exp(s x_\al)  \;\,,\,\;
   \phi_\al\left(\begin{array}{cc}
  1 & 0\\
  s & 1
  \end{array}\right) \;:=\; \exp(sx_{-\al})\;\,,\,\;(s\in \F^\times)\;.
\end{eqnarray*}
{\bf (b)} Denote by $N$ the subgroup generated by $T$ and $
  n_\al := \phi_\al\left(\begin{array}{cc}
        0      & 1  \\
  -1 & 0
  \end{array}\right)$, $\al\in \Delta_{re}$. The Weyl group $\We$ can be identified with the group $N/T$ by the isomorphism 
$\kappa:\,\We \to  N/T$ given by $\kappa(\sigma_\al):=n_\al T$, $\al\in \rW $.
We denote an arbitrary element $n\in N$ with $\kappa^{-1}(nT)=\sigma\in\We$  by $n_\sigma$. The set of weights $P(V)$ of an
admissible $\g$-module $(V,\pi)$ is $\We$-invariant, and $ n_\sigma V_\la =V_{\sigma \la}$, $\la\in P(V)$.
\vspace*{1ex}\\
{\bf (c)} Let $U^\pm$ be the subgroups generated by $U_\al$, $\al\in\Delta_{re}^\pm $. Let $U_f:=\exp(\n_f)$. Then $U^\pm$ and $U_f$ are
normalized by $T$. Set
\begin{eqnarray*}
   B^\pm\;\,:=\,\; T\ltimes U^\pm &\quad,\quad &  B_f\;\,:=\,\; T\ltimes U_f\;\;.
\end{eqnarray*}
The pairs ($B^\pm$, $N$) are  twinned BN-pairs of $G$ with the property $B^+\cap B^-= B^\pm\cap N = T$. The pair 
($B_f$, $N$) is a BN-pair of $G_f$ with $B_f\cap N=T$. There are the {\it Bruhat} and 
{\it Birkhoff decompositions}
\begin{eqnarray*}
  G \;\,=\;\, \dot{ \bigcup_{\sigma\in{\cal W}}} B^\epsilon\sigma B^\delta \quad,\quad 
  G_f \;\,=\;\, \dot{ \bigcup_{\sigma\in{\cal W}}} B^\epsilon\sigma B_f \qquad,
  \qquad \epsilon,\delta\;\in\;\{\,+\,,\,-\,\}\;\;.
\end{eqnarray*}
{\bf (d)} There are also Levi decompositions of the standard parabolic subgroups. In this article we only use the corresponding
decompositions for the groups $U^\pm$ and $U_f$:
Set $\W^\pm_J:=\W^\pm\cap \sum_{j\in J}\Z\,\al_j$, and ${(\W^J)}^\pm:=\W^\pm\setminus\sum_{j\in J}\Z\,\al_j$. Similarly define 
$(\W_J)^\pm_{re}$ and $(\W^J)^\pm_{re}$ by replacing $\W^\pm$ by $\W_{re}^\pm$. Set 
$(\n_J)^\pm:=\bigoplus_{\al\in \Delta_J^\pm}\g_\al$, $(\n_f)_J:=\prod_{\al\in\Delta_J^+}\g_\al$, and 
$(\n_f)^J:=\prod_{\al\in(\Delta^J)^+}\g_\al$. We have
\begin{eqnarray*}
   U^\pm\;\,=\;\, U^\pm_J \ltimes (U^J)^\pm  &\quad,\quad &  U_f\;\,=\;\, (U_f)_J \ltimes (U_f)^J\;\;.
\end{eqnarray*}
Here $U_J^\pm$ is the group generated by $U_\al$, $\al\in (\W_J)^\pm_{re}$. $(U^J)^\pm$ is the smallest normal subgroup of $U^\pm$
containing $U_\al$, $\al\in (\W^J)^\pm_{re}$. This group equals $\,\bigcap_{\sigma\in{\cal W}_J}\sigma U^\pm\sigma^{-1}$.
Furthermore $(U_f)_J:=\exp((\n_f)_J)$ and $(U_f)^J:=\exp((\n_f)^J)$.\vspace*{1ex}\\ 
For $w\in\We$ set  $U_w := U\cap w U^- w^{-1}$, $U^w:=U\cap  w U w^{-1}$, and $(U_f)^w :=U_f\cap  w U_f w^{-1}$.
Then $U_w=\prod_{\al\in \Phi_w }U_\al$ (arbitrary order of the product), where $\Phi_w= \prW\cap w \nrW= 
\pW\cap w\nW$. The 
multiplication maps $U_w\times U^w\to U$ and $U_w\times (U_f)^w\to U_f$ are bijective. Set $U_w^-:=(U_w)^*$ 
and $(U^w)^-:=(U^w)^*$.\vspace*{1ex}\\
The derived minimal Kac-Moody group $G'$ is identical with the Kac-Moody group as defined in \cite{KP1}. It is generated by 
the root groups $U_\al$, $\al\in \rW$. We have $G = G'\rtimes T_{rest}$, where $T_{rest}:=H_{rest}\otimes_{\mathbb Z} \F$ is a 
subtorus of $T$.
The group $G_f$ is identical with the Kac-Moody group of \cite{Sl1} for a simply connected minimal free realization.
\vspace*{1ex}\\
The monoid $\GD$ has the following important structural properties:
The Kac-Moody group $G$ is the unit group of $\GD$. Every idempotent is $G$-conjugate to some idempotent $e(R(\Th))$, $\Th$ special. We have
\begin{eqnarray*}
  \GD&=& \dot{\bigcup_{\Th\,special}} G e(R(\Th)) G\;\;.
\end{eqnarray*}
We get an abelian submonoid of $\GD$ by $\TD:=\dot{\bigcup}_{R\in{\cal R}(X)} T e(R) $. We get a submonoid of $\GD$ by 
$\ND:=\dot{\bigcup}_{R\in{\cal R}(X)} N e(R)$. Define a congruence relation on $\ND$ as follows:
\begin{eqnarray*}
\quad\hat{n}\:\sim\:\hat{n}'\quad:\iff\quad\hat{n}T\:=\:\hat{n}'T\quad\iff \quad\hat{n}'\:\in\:\hat{n}T\quad\iff\quad\hat{n}\in\hat{n}'T\;\;.
\end{eqnarray*}
The Weyl monoid $\WeD$ is isomorphic to the monoid $\ND/T$, an isomorphism $\kappa:\WeD\to \ND/T$ given by $\kappa(\sigma\ve{R})=n_\sigma e(R) T$.
\vspace*{1ex}\\
$\GD$ has {\it Bruhat} and {\it Birkhoff decompositions}: 
\begin{eqnarray*}
  \GD &=& \dot{ \bigcup_{\hat{n}\in \widehat{N}}}\,\; U^\epsilon\,\hat{n}\, U^\delta \;\,=\;\,
   \dot{ \bigcup_{\hat{\sigma}\in \widehat{\cal W}}}\,\; B^\epsilon\,\hat{\sigma}\, B^\delta \quad,\quad  \epsilon,\,\delta \;\in \;\{\,+\,,\,-\,\}\;\;.
\end{eqnarray*}
For later reference we state the following formulas, which are useful for computations in $\GD$:\\
{\bf ($\al$)} Let $R$, $S$ be faces of the Tits cone, and $n_\sigma\in N$. Then
\begin{eqnarray*}
 e(R)e(S)=e(R\cap S) \qquad ,\qquad  n_\sigma e(R) n_\sigma^{-1}\;\,=\;\, e(\sigma R)\;\;.
\end{eqnarray*}  
{\bf ($\beta$)} An element $g$ of $T$, $N$, $U$, $U^-$, resp. $G$ satisfies 
\begin{eqnarray*}
   e(R(\Th)) g &=& e(R(\Th))
\end{eqnarray*}
if and only if it satisfies
\begin{eqnarray*}
    g^* e(R(\Th))  &=& e(R(\Th))
\end{eqnarray*}
if and only if it is contained in $T_\Th$, $N_\Th$, $U_\Th$, $U_\Th^-\ltimes (U^{\Th\cup\Th^\bot})^{-}$, resp. 
$G_\Th\ltimes (U^{\Th\cup\Th^\bot})^-$.
Here $T_\Th$ is the subtorus of $T$ generated by $t_{h_j}(s)$, $j\in \Th$, $s\in\F^\times$, $N_\Th$ is the subgroup of $N$ generated by $T_\Th$ and 
$n_{\al_j}$, $j\in \Th$, and $G_\Th$ is the subgroup of $G$ generated by $U_{\al_j}^\pm$, $j\in \Th$.\vspace*{1ex}\\
{\bf ($\gamma$)} An element $g$ of $T$, $N$, $U$, $U^-$, resp. $G$ satisfies 
\begin{eqnarray*}
   g e(R(\Th)) g^{-1} &=& e(R(\Th))
\end{eqnarray*}
if and only if it is contained in the groups $T$, $ N_{\Th\cup\Th^\bot} T$,  $U_{\Th\cup\Th^\bot}$, $U_{\Th\cup\Th^\bot}^-$, resp. 
$G_{\Th\cup\Th^\bot} T$.\vspace*{1ex}\\ 
{\bf ($\delta$)} In particular we have
\begin{eqnarray*}
    U e(R(\Th))  &=& U_{\Th^\bot} e(R(\Th))\;\,=\;\, 
    e(R(\Th))U_{\Th^\bot}\;\;,\\
     e(R(\Th))U^-  &=& e(R(\Th))U^-_{\Th^\bot}\;\,=\;\, 
      U^-_{\Th^\bot}e(R(\Th))\;\;.
\end{eqnarray*}
%
%
%
%
%
%
%
%
{\bf Sets with coordinate rings:}
%
%
%
%
%
We call a point separating algebra of functions $\FK{A}$ on a set $A$ a {\it coordinate ring}.
The closed sets of the {\it Zariski topology} on $A$ are given by the zero sets of the functions of $\FK{A}$. 
The set $A$ is irreducible if and only if $\FK{A}$ is an integral domain.\\
A morphism of sets with coordinate rings ($A,\FK{A}$) and ($B,\FK{B}$) consists of a map $\phi:A\to B$, whose comorphism $\phi^*:\FK{B}\to\FK{A}$ 
exists. In particular a morphism is Zariski continuous.\vspace*{1ex}\\
If ($B,\FK{B}$) is a set with coordinate ring, and $A$ is a nonempty subset of $B$, we get a coordinate ring on 
$A$ by restricting the functions of $\FK{B}$ to $A$.\\
If ($A,\FK{A}$) is a set with coordinate ring and $f\in\FK{A}\setminus\{0\}$, the principal open set 
$D_A(f):=\Mklz{a\in A}{f(a)\neq 0}$ is equipped with a coordinate ring by identifying the localization $\FK{A}_f$ in 
the obvious way with an algebra of functions on $D_A(f)$. The principal open set $D_A(f)$ is irreducible if and 
only if $A$ is irreducible.\\
If ($A,\FK{A}$) and ($B,\FK{B}$) are sets with coordinate rings, then the product $A\times B$ is equipped with a coordinate ring by identifying the tensor product 
$\FK{A}\otimes\FK{B}$ in the obvious way with an algebra of functions on $A \times B$. The product $A\times B$ 
is irreducible if and only if $A$ and $B$ are irreducible.\vspace*{1ex}\\
%
%
%
{\bf The coordinate ring of $\GD$:}
%
%
By the Tannaka Krein reconstruction given in \cite{M1}, Section 4, the monoid $\GD$ is equipped with a natural coordinate ring. It can also be defined as follows:
For a module $V$ of ${\cal O}_{adm}$, $v,w\in V$, and $\kBl$ a nondegenerate symmetric contravariant bilinear form on $V$, call the function 
$f_{vw}:\,\GD\to\F$ defined by $f_{vw}(x):=\kB{v}{xw}$, $x\in\GD$, a {\it matrix coefficient} of $\GD$. The set of all such matrix 
coefficients $\FK{\GD}$ is a coordinate ring on $\GD$, which is an integral domain. For a set $M\subseteq \GD$ we denote by $\overline{M}$ its Zariski closure.\\ 
The Chevalley involution $*:\GD\to\GD$ is a morphism. We denote its comorphism also by $*:\FK{\GD}\to\FK{\GD}$ and 
call it {\it Chevalley involution}.
Right and left multiplications with elements of $\GD$ are morphisms. In particular an action of 
$\GD\times\GD$ on $\GD$ from the right by morphisms is given by
\begin{eqnarray*}
 x \,(g,h):= g^* x h && x,g,h\in\GD\;\;.
\end{eqnarray*}
The comorphisms induce an action of $\GD\times\GD$ on $\FK{\GD}$ from the left.\\
In this article we fix a nondegenerate symmetric contravariant bilinear form on $L(\La)$ for every $\La\in P^+$.
The coordinate ring $\FK{\GD}$ admits a {\it Peter-Weyl theorem}: The map $\bigoplus_{\La\in P^+} L(\La)\otimes L(\La) \to  \FK{\GD}$ induced by 
$v\otimes w\mapsto f_{vw}$ is an isomorphism of $\GD\times \GD$-modules. It identifies the direct sum of the switch 
maps of the factors with the Chevalley involution.\vspace*{1ex}\\
The algebra of {\it strongly regular functions} $\FK{G}$ is obtained by restricting the functions of $\FK{\GD}$ onto 
$G$. The restriction map is an isomorphism from $\FK{\GD}$ to $\FK{G}$.
Restricting the functions of $\FK{G}$ onto $G'$, resp. $T_{rest}$ gives the algebras $\FK{G'}$, resp. $\FK{T_{rest}}$ 
the first identical with the algebra of strongly regular functions as defined in \cite{KP2}, the second the classical 
coordinate ring of the torus $T_{rest}\,$. The multiplication map $G'\times T_{rest}\to G$ is an isomorphism.
\vspace*{1ex}\\   
The monoids $\TD$, $\ND$, $\GD$ are the Zariski closures of $T$, $N$, $G$, and $G$ is the Zariski open dense unit group of $\GD$.\\
Denote by $\leq$ the Bruhat order on $\We$. Lemma 3.4 of \cite{KP2} is also valid for the slightly enlarged group $G$, which we use here. It gives 
the relative closures of Bruhat and Birkhoff cells of $G$:  
\begin{eqnarray}\label{clBBG}
   \overline{B^\eps w B^\eps}\cap G \;\,=\;\, \dot{\bigcup_{w'\in {\cal W}\atop w'\leq w}} B^\eps w' B^\eps
  &\mb{ and }& \overline{B^- w B} \cap G\;\,=\;\, \dot{\bigcup_{ w'\in {\cal W}\atop w'\geq  w}} B^- w' B\;\;,
\end{eqnarray}
where $\eps\in\{+,-\}$ and $w\in\We$.\vspace*{1ex}\\
%
%
%
%
{\bf The spectrum of $\F$-valued points of $\FK{\GD}$:}
We denote the {\it $\F$-valued points} of $\FK{\GD}$, i.e., the homomorphisms of algebras from $\FK{\GD}$ to $\F$, by $\Spm\FK{\GD}$.\\
A function $f\in\FK{\GD}$ induces a function on $\Spm\FK{\GD}$, assigning to $\phi\in\Spm\FK{\GD}$ the value $\phi(f)$. We denote this function also 
by $f$. In this way $\Spm\FK{\GD}$ is equipped with a coordinate ring isomorphic to $\FK{\GD}$. Its Zariski topology coincides with the relative 
topology induced by the topology of the spectrum of $\FK{\GD}$.
For a set $M\subseteq \Spm\FK{\GD}$ we denote by $\overline{M}^{spm}$ its Zariski closure.\vspace*{1ex}\\
As a set, the $\F$-valued points of $\FK{\GD}$ can be described as follows: There is a surjective map 
$\tr: \GfD\times\GfD \to \Spm\FK{\GD}$ given by 
\begin{eqnarray*}
   (x\tr y)(f_{vw}) \;\,:=\;\, \kB{xv}{yw} &,& x,y\in\GfD,\;\;v,w\in L(\La),\;\;\La\in P^+\;\;.
\end{eqnarray*}
The set of fibres of this map coincides with the partition of $\GfD\times\GfD$ corresponding to the equivalence 
relation, which is generated by 
\begin{eqnarray*} 
   (x,\,zy) \;\sim\; (z^*x,\,y) &,& x,y\in\GfD,\;\;z\in\GD\;\;.
\end{eqnarray*}
In particular $x\tr zy=x z^*\tr y$, $x,y\in\GfD$, $z\in\GD$.\\
The monoid $\GfD\times\GfD$ acts in a natural way on $\Spm\FK{\GD}$ by morphisms from the right. The map 
$\tr$ is equivariant, i.e.,  
\begin{eqnarray*}
   (x\tr y)(\ti{x},\ti{y}) \;\,=\;\, x\ti{x}\tr y\ti{y} &,& x,\ti{x}, y,\ti{y}\in\GfD\;\;.
\end{eqnarray*}
The Chevalley involution of $\FK{\GD}$ induces an involutive morphism $*$ on $\Spm\FK{\GD}$, also called {\it Chevalley involution}. Using 
the map $\tr$ it can be described by
\begin{eqnarray*}
   (x\tr y)^*\;\,=\;\, y \tr x &,& x,y\in\GfD\;\;.
\end{eqnarray*}
In \cite{M2} we investigated the $G_f\times G_f$-orbit decomposition (closure relation of the orbits, irreducibility of the orbits, coverings of the 
orbits by big cells, transversal stratified slices to the orbits). As an aid for these investigations, we showed the {\it Birkhoff decomposition}
\begin{eqnarray*}
 \Spm\FK{\GD}&=& \dot{\bigcup_{\hat{w}\in{\cal W}}} B_f\tr \hat{w}B_f\;\;.
\end{eqnarray*}
Note that $\GfD$ embeds into $\Spm\FK{\GD}$ by assigning $x\in\GfD$ the element $1\tr x$. In this way this decomposition of $\Spm\FK{\GD}$ extends one of the Birkhoff decompositions of $\GD$, and a corresponding Birkhoff decomposition, 
which holds for $\GfD$.\vspace*{1ex}\\
%
%
%
%
{\bf The Kostant cones:}
Fix $\La\in P^+$. Recall that we have fixed a nondegenerate symmetric contravariant bilinear form $\kBl$ on 
$L(\La)$. Equip $L(\La)$ with the coordinate ring $\FK{L(\La)}$ generated by the matrix coefficients $f_v$, 
$v\in L(\La)$, where $f_v(w):=\kB{v}{w}$ for all $w\in L(\La)$. It is a symmetric algebra in the linear space given by these functions.
For a set $M\subseteq L(\La)$ denote by $\overline{M}$ its Zariski closure.\vspace*{1ex}\\ 
The set ${\cal V}_\La:=G (L(\La)_\La)\subseteq L(\La)$ is Zariski closed. It is called {\it Kostant cone}.\vspace*{1ex}\\
Now assume that $\La\in F_J$, $J\subseteq I$. Denote by $\We^J$ the set of minimal coset representatives of 
$\We/\We_J$. Because the parabolic subgroup $\We_J$ is the 
$\We$-stabilizer of $\La$, the evaluation map $\We^J \to \We^J\La=\We\La$ is bijective. By this map the Bruhat order on $\We^J$ induces an order on $\We\La$. The corresponding inverse order on $\We\La$ is denoted by $\preceq $. 
If $\la,\,\mu\in\We\La$ then $\la\preceq\mu$ implies $\la\leq \mu$.\vspace*{1ex}\\
For $v\in L(\La)$ denote by $supp(v)$ the set of weights of the nonzero weight space components of $v$. Denote by 
$S(v)$ the convex hull of $supp(v)$ in $\h_{\mathbb R}^*$.\\
For $v\in{\cal V}_\La\setminus\{0\}$ the vertices of $S(v)$ are given by $S(v)\cap\We\La$. The edges of $S(v)$ are parallel to real roots. The two vertices of an edge are comparable in $\preceq$. S(v) has one maximal and 
one minimal vertex. \vspace*{1ex}\\
For $\la\in \We\La$ set
\begin{eqnarray*}
   {\cal V}_\La^+(\la) &:=& \Mklz{v\in {\cal V}_\La\setminus\{0\} }{\la \mb{ is the minimal vertex of } S(v)}\;\;,\\
   {\cal V}_\La^-(\la) &:=& \Mklz{v\in {\cal V}_\La\setminus\{0\} }{\la \mb{ is the maximal vertex of } S(v)}\;\;.
\end{eqnarray*}
Then ${\cal V}_\La(\la)^\eps=U^\eps (L(\La)_\la\setminus\{0\})$, and 
${\cal V}_\La\setminus\{0\}=\dot{\bigcup}_{\la\in{\cal W}\La}{\cal V}_\La(\la)^\eps$ where $\eps\in\{+,-\}$, and
\begin{eqnarray*}
  \overline{{\cal V}_\La(\la)^+}\setminus\{0\} \;=\; \dot{\bigcup_{\mu\in{\cal W}\La,\,\mu\succeq \la }} {\cal V}_\La(\mu)^+ \quad,\quad
  \overline{{\cal V}_\La(\la)^-}\setminus\{0\} \;=\; \dot{\bigcup_{\mu\in{\cal W}\La,\,\mu\preceq \la }} {\cal V}_\La(\mu)^- \;\;.
\end{eqnarray*}
The Kostant cones are the affine cones of the flag varieties. For our investigations they are more important than 
the flag varieties, because the monoid $\GD$ acts on the Kostant cones by morphisms, but it does not act on the 
flag varieties. In particular ${\cal V}_\La=\GD (L(\La)_\La)$.  
%
%
%
%
\section{Extensions of the Bruhat order\label{SBO}}
%
%
%
Since $\GD\times\GD$ acts on $\GD$ by morphisms, the closures of the $B^{-\eps}\times B^\de$-orbits, $\eps,\,\de\in\{+,-\}$, are unions of 
$B^{-\eps}\times B^\de$-orbits. Because of the Bruhat and Birkhoff decompositions
\begin{eqnarray*}
  \GD &=& \dot{\bigcup_{\hat{w}\in \widehat{\cal W}}}\, B^\eps\hat{w} B^\de
\end{eqnarray*}
the closure relations of the $B^{-\eps}\times B^\de$-orbits determine relations on $\WeD$:
\begin{Def}\label{BO1} For $\hat{w},\hat{w}'\in\WeD$ and $\eps\in \{+,-\}$ define:
\begin{eqnarray*}
   \hat{w}' \leq_{\eps\eps} \hat{w} &:\iff&  B^\eps\hat{w}' B^\eps\;\subseteq\;\overline{B^\eps\hat{w} B^\eps}\;.\\ 
   \hat{w} \leq_{-+} \hat{w}' &:\iff&  B^-\hat{w}' B\;\subseteq\;\overline{B^-\hat{w} B}\;.
\end{eqnarray*}
\end{Def}
{\bf Remarks:}
{\bf (1)} Due to this definition
\begin{eqnarray*}
  \overline{B^\eps\hat{w} B^\eps} \;\,=\;\, \dot{\bigcup_{\hat{w}'\in\widehat{\cal W}\atop \hat{w}'\leq_{\eps\eps} \hat{w}}} B^\eps\hat{w}' B^\eps
  &\mb{ and }& \overline{B^-\hat{w} B} \;\,=\;\, \dot{\bigcup_{\hat{w}'\in\widehat{\cal W}\atop \hat{w}'\geq_{-+} \hat{w}}} B^-\hat{w}' B\;.
\end{eqnarray*}
{\bf (2)} Due to equations (\ref{clBBG}) the relations $\leq_{++}$, $\leq_{--}$, and $\leq_{-+}$ extend the Bruhat order $\leq$ of the Weyl group $\We$.\\
The relations $\leq_{++}$, $\leq_{--}$, and $\leq_{-+}$ are the closure relations of orbit decompositions. Therefore these relations are 
reflexive and transitive. But we can not conclude that they are antisymmetric, because for this we would need to know that the orbits are locally 
closed. \\ 
{\bf (3)} The Chevalley involution $*:\GD\to\GD$ is a morphism with $(B^\eps \hat{w} B^\delta)^*=B^{-\delta} \hat{w}^{inv} B^{-\eps}$, $\hat{w}\in\WeD$, 
$\eps,\delta\in\{+,-\}$. From the definition of the extended Bruhat orders follows immediately, that the inverse map $\mb{}^{inv}:\We\to\We$ is an 
isomorphism of the relations $(\We,\leq_{\eps\eps})$ and  $(\We,\leq_{-\eps -\eps})$, $\eps=+,-$. It is an automorphism of $(\We,\leq_{-+})$.\\
{\bf (4)} In this article we do not treat the Birkhoff cells $B\hat{w}B^-$, $\hat{w}\in\WeD$. These cells can not be treated in a similar way as the Birkhoff cells $B^-\hat{w}B$, $\hat{w}\in\WeD$, because the coordinate ring 
$\FK{\GD}$ does not contain the matrix coefficients of the admissible lowest weight representations. Even the 
relative closures of $BwB^-$, $w\in\We$, in the Kac-Moody group $G$ have not been determined.\vspace*{1ex}\\
Our first aim is to determine these relations explicitely. To this end we introduce three normal forms of the elements of $\WeD$, which are similar to the standard form of an element of a Renner monoid introduced in \cite{PePuRe}.\\
For $J\subseteq I$ denote by $\We^J$ the minimal coset representatives of $\We/\We_J$, and denote by $\mb{}^J\We$ the minimal coset 
representatives of $\We_J\backslash \We$. 
\begin{Prop}\label{BO2} Let $\hat{w}\in\WeD$.\\
(I) There exists a uniquely determined special set $\Th$, and uniquely determined elements $w_1\in\We^\Th$, $w_2\in \mb{}^{\Th\cup\Th^\bot}\We$, 
such that
\begin{eqnarray*}
  \hat{w}  & = &  w_1\ve{R(\Th)} w_2\;\;.
\end{eqnarray*}
(II) There exists a uniquely determined special set $\Th$, and uniquely determined elements $w_1\in\We^{\Th\cup\Th^\bot}$, $w_2\in \mb{}^\Th \We$, 
such that
\begin{eqnarray*}
  \hat{w}  & = &  w_1\ve{R(\Th)} w_2\;\;.
\end{eqnarray*}
(III) There exists a uniquely determined special set $\Th$, and uniquely determined elements $w_1\in\We^{\Th\cup\Th^\bot}$, $w_2\in\We_{\Th^\bot}$, 
$w_3\in \mb{}^{\Th\cup\Th^\bot} \We$, such that
\begin{eqnarray*}
  \hat{w}  & = &  w_1 w_2\ve{R(\Th)} w_3\;\,=\;\, w_1 \ve{R(\Th)}w_2 w_3 \;\;.
\end{eqnarray*}
\end{Prop}
{\bf Remarks:} 
{\bf (1)} By applying the inverse map $\mb{}^{inv}:\WeD\to\WeD$ to an element $\hat{w}$ in normal form I resp. II 
we obtain the element $\hat{w}^{inv}$ in 
normal form II resp. I. By applying this map to an element $\hat{w}$ in normal form III we obtain the element $\hat{w}^{inv}$ in normal form III.\\ 
{\bf (2)} As shown in the following proof, the multiplication map of $\We$ restricts to bijective maps  
\begin{eqnarray}\label{multW}
    \We^{\Th\cup\Th^\bot}\times \We_{\Th^\bot}\;\to\; \We^\Th \qquad \mb{and}\qquad 
    \We_{\Th^\bot}\times \mb{}^{\Th\cup\Th^\bot}\We\;\to\;\mb{}^\Th\We\;\;.
\end{eqnarray} 
Therefore we can read off the normal forms I and II from the normal form III.\vspace*{1ex}\\ 
\Proof {\bf 1)} We first show the existence and uniqueness of normal form I: We have $\WeD = 
\dot{\bigcup}_{\Th\;special} \We \ve{R(\Th)}\We$, and by using formulas (\ref{WeDF1}) and (\ref{WeDF2}) we find
\begin{eqnarray*}
        \We \ve{R(\Th)}\We \;\,=\;\, \We^{\Th\cup\Th^\bot}\We_{\Th\cup\Th^\bot} \ve{R(\Th)}\We \;\,=\;\, 
        \We^{\Th\cup\Th^\bot}\ve{R(\Th)}\We_{\Th\cup\Th^\bot} \We \\ \;\,=\;\, \We^{\Th\cup\Th^\bot}\ve{R(\Th)}\We \;\,=\;\, 
        \We^{\Th\cup\Th^\bot}\ve{R(\Th)}\We_\Th \We^\Th \;\,=\;\, \We^{\Th\cup\Th^\bot}\ve{R(\Th)}\We^\Th \;\;.
\end{eqnarray*}
To show the uniqueness let $w_1\ve{R(\Th)}w_2$ and $w_1'\ve{R(\Th')}w_2'$ be normal forms I of the same element of 
$\WeD$. Then by using equation (\ref{wRwinv}) we find
\begin{eqnarray*}
   \ve{w_1 R(\Th)}w_1 w_2 \;\,=\;\, w_1\ve{R(\Th)}w_2\;\,=\;\, w_1'\ve{R(\Th')}w_2' \;\,=\;\, 
   \ve{w_1' R(\Th')}w_1' w_2'\;\;, 
\end{eqnarray*} 
which is equivalent to 
\begin{eqnarray*}
  w_1 R(\Th)= w_1' R(\Th') \qquad\mb{and}\qquad  w_1' w_2' (w_1 w_2)^{-1}\in w_1'\We_{\Th'}(w_1')^{-1}\;\;.
\end{eqnarray*} 
From the first equation follows 
$\Th=\Th'$ and $w_1\We_{\Th\cup\Th^\bot}=w_1'\We_{\Th\cup\Th^\bot}$. Since the minimal coset representatives $w_1$, 
$w_1'$ are uniquely determined we find $w_1=w_1'$.
Inserting in the second expression we get $w_2'w_2^{-1}\in \We_{\Th}$, resp. $\We_\Th w_2'=\We_\Th w_2$. Because the 
minimal coset representatives $w_2$, $w_2'$ are uniquely determined this implies $w_2=w_2'$.\vspace*{1ex}\\ 
{\bf 2)} The existence and uniqueness of normal form II follows from the existence and uniqueness of normal form I by using the inverse map $\mb{}^{inv}:\WeD\to\WeD$.\vspace*{1ex}\\ 
{\bf 3)} If we show the bijectivity of the restricted multiplication maps (\ref{multW}), then the existence and 
uniqueness of normal form I and II together with formula (\ref{WeDF1}) imply the existence and uniqueness of normal form III.
We only have to show the bijectivity of the first map, the bijectivity of the second follows by applying the 
inverse map of $\We$.\\
We have $\We=\We^{\Th\cup\Th^\bot}\We_{\Th\cup\Th^\bot}=\We^{\Th\cup\Th^\bot}\We_{\Th^\bot}\We_\Th$, and the 
corresponding multiplicative decomposition of the elements are unique. It remains to show that 
$\We^{\Th\cup\Th^\bot}\We_\Th\subseteq \We^\Th$.
If $w_1\in \We^{\Th\cup\Th^\bot}$ and $w_2\in \We_{\Th^\bot}$ then for all $i\in \Th$ we have 
\begin{eqnarray*}
  w_1 w_2\al_i\;\,=\,\; w_1\al_i\;\,\in\;\,\prW\;\;.
\end{eqnarray*}
Therefore $ w_1 w_2\in \We^\Th$.\\
\End
The following theorem describes the relations $\leq_{++}$, $\leq_{--}$, and $\leq_{-+}$ explicitely. The 
description of 1a) (iii) is similar to the result for reductive algebraic monoids obtained in \cite{PePuRe}.\\ 
For $J\subseteq I$ and $w\in \We$ denote by $w^J$ the minimal coset representative of $w\We_J$, and by $\mb{}^J w$ the minimal coset 
representative of $\We_J w$. 
\begin{Theorem}\mb{}\label{BO3}\\
{\bf 1a)} Let $\hat{w}= w_1 \ve{R(\Th)} w_2 $, $\hat{w}' =w'_1 \ve{R(\Th')} w'_2$ be elements of $\WeD$ in normal form I. Then the following 
statements are equivalent:\\ 
(i) $\hat{w}\leq_{++}\hat{w}' $\\
(ii) $\Th\supseteq\Th'$ and there exists an element $w\in\We_{{\Th'}^\bot}$ such that 
$w_1\leq (w_1' w)^\Th$ and $w_2\geq \mb{}^\Th(w^{-1}w_2')$.\\
(iii) $\Th\supseteq\Th'$ and there exists an element $w\in\We_{{\Th'}^\bot}\We_\Th$ such that 
$w_1\leq w_1' w $ and $w_2\geq w^{-1}w_2'$.\vspace*{1ex}\\
{\bf b)} Let $\hat{w}= w_1 \ve{R(\Th)} w_2 $, $\hat{w}' =w'_1 \ve{R(\Th')} w'_2$ be elements of $\WeD$ in normal form II. Then the following 
statements are equivalent:\\ 
(i) $\hat{w}\leq_{--}\hat{w}' $\\
(ii) $\Th\supseteq\Th'$ and there exists an element $w\in\We_{{\Th'}^\bot}$ such that 
$w_1\geq (w_1' w)^\Th$ and $w_2\leq \mb{}^\Th(w^{-1}w_2')$.\\
(iii) $\Th\supseteq\Th'$ and there exists an element $w\in\We_{{\Th'}^\bot}\We_\Th$ such that 
$w_1\geq w_1' w$ and $w_2\leq w^{-1}w_2'$.\vspace*{1ex}\\
{\bf 2)} Let $\hat{w}= w_1 \ve{R(\Th)} w_2 $, $\hat{w}' =w'_1 \ve{R(\Th')} w'_2$ be elements of $\WeD$ both in normal form I, or both in normal 
form II. Then the following statements are equivalent:\\
(i) $\hat{w}\geq_{-+}\hat{w}' $ \\
(ii) $\Th\supseteq\Th'$ and there exists an element $w\in\We_{{\Th'}^\bot}$ such that $w_1\geq (w_1' w)^\Th$ and 
$w_2\geq \mb{}^\Th(w^{-1}w_2')$.
\end{Theorem}
\Proof
We only have to prove the statements of the theorem, which use normal form I. Then the statements which use 
normal form II follow by applying the inverse map $\mb{}^{inv}:\WeD\to\WeD$, compare Remark (3) following 
Definition \ref{BO1}, and Remark (1) following Proposition \ref{BO2}.\vspace*{1ex}\\
{\bf (a)} First we show the equivalence of (ii) and (iii) of 1a): 
Let (ii) be valid. Choose an element $\ti{w}\in\We_\Th$ such that $\mb{}^\Th(w^{-1}w_2')=\ti{w}^{-1}w^{-1}w_2'$. Then we have
\begin{eqnarray*}
  w_1\leq (w_1'w)^\Th = (w_1'w\ti{w})^\Th\leq w_1'w\ti{w} \;,\\
  w_2\geq \mb{}^\Th(w^{-1}w_2') = (w\ti{w})^{-1}w_2'\;.
\end{eqnarray*}
Let (iii) be valid. Write the element $w$ in the form $w=\ti{w}\ti{w}'$ with $\ti{w}\in\We_{\Th'^\bot}$ 
and $\ti{w}'\in \We_\Th$. From $w_1\leq w_1' w$ follows $w_1=w_1^\Th\leq (w_1'w)^\Th=(w_1'\ti{w})^\Th$. Similarly we get 
$w_2\geq \mb{}^\Th(\ti{w}^{-1}w_2')$.\vspace*{1ex}\\
The easy part in the following proof of the equivalence (i) and (ii) in 1), 2a) is the direction from (ii) to (i). 
The proof of this direction, which uses only the formulas for the relative closures of the Bruhat and Birkhoff 
cells of a Kac-Moody group, is for 1a) similar, and for 2) a modification of the corresponding proof in 
\cite{PePuRe}. The proof of the opposite direction is quite different.\vspace*{1ex}\\ 
{\bf (b)} To prepare the proof of the direction from (i) to (ii) of 1a) and 2) first note: 
Let $\La\in P^+$, and let $v,v'\in {\cal V}_\La$ such that $\kB{v}{v'}\neq 0$. 
Since different weight spaces of $L(\La)$ are $\kBl$-orthogonal we have $supp(v)\cap supp(v')\neq\emptyset$. Let 
$\la$ be an element of this intersection. Then the biggest vertex of $S(v)$ is bigger than $\la$, which is 
bigger than the smallest vertex of $S(v')$ (with respect to $\preceq$).
\vspace*{1ex}\\
{\bf (c)} Now we can prove the direction from (i) to (ii) of 1a) and 2): Let $\eps\in\{+,-\}$ 
and $B^\eps\hat{w}B\subseteq \overline{B^\eps\hat{w}'B}$.\\ 
Choose an element $\La\in F_\Th\cap P$. Choose elements $v_{w_1\La}\in  L(\La)_{w_1\La}\setminus\{0\}$ and 
$v_{w_2\La}\in  L(\La)_{w_2\La}\setminus\{0\}$. To cut short the notation set 
$g:=f_{v_{w_1\La} v_{(w_2)^{-1}\La}}$.\\
Since $\La\in R(\Th)$, for any element $\hat{n}_{\hat{w}}\in \hat{N}$ belonging to $\hat{w}$ we have
\begin{eqnarray*}
   g(\hat{n}_{\hat{w}}) \;\,=\;\, \kB{v_{w_1\La}}{\hat{n}_{w_1\varepsilon(R(\Th))w_2} v_{(w_2)^{-1}\La}}\;\,\neq\;\,0 \;\;.
\end{eqnarray*}
Therefore $g$ does not vanish entirely on the closure $\overline{B^\eps\hat{w}'B}$, which implies that it also can not vanish entirely on $B^\eps\hat{w}'B$.
Therefore there exists an element $\hat{n}_{\hat{w}'}\in \hat{N}$ belonging to $\hat{w}'$, and elements $u_\eps\in U^\eps$, $u\in U$, such that
\begin{eqnarray*}
   g(u_\eps \hat{n}_{\hat{w}'}u) \;\,=\;\, \kB{(u_\eps)^* v_{w_1\La}}{\hat{n}_{w_1'\varepsilon(R(\Th'))w_2'} u v_{(w_2)^{-1}\La}}\;\,\neq\;\,0 \;\;.
\end{eqnarray*}
The vertices of $S(u v_{(w_2)^{-1}\La})$ are of the form
\begin{eqnarray*}
  \ti{w}\La \quad\mb{with}\quad \ti{w}\in\We\quad\mb{such that}\quad \ti{w}^\Th\leq ((w_2)^{-1})^\Th=(w_2)^{-1}\;\;.
\end{eqnarray*} 
Therefore the vertices of $S(\hat{n}_{w_1'\varepsilon(R(\Th'))w_2'}u v_{(w_2)^{-1}\La})$ are of the form
\begin{eqnarray*}
  w_1' w_2'\ti{w}\La \quad\mb{with}\quad \ti{w}\in\We\quad\mb{such that}\quad \ti{w}^\Th\leq (w_2)^{-1}
  \quad\mb{and}\quad w_2'\ti{w}\La\in R(\Th')\;\;.
\end{eqnarray*} 
Here the inequality $\ti{w}^\Th\leq (w_2)^{-1}$ is equivalent to $\mb{}^\Th(\ti{w}^{-1})\leq w_2$. By comparing the type of facets in 
the formula $w_2'\ti{w}\La\in R(\Th')$ on the left and on the right we find 
\begin{eqnarray*}
   F_\Th\subseteq \overline{F_{\Th'}} \quad\mb{resp.}\quad \Th\supseteq\Th'\quad\mb{resp.}\quad R(\Th)\subseteq R(\Th')\;\;.
\end{eqnarray*}
We also get $w_2'\ti{w}\in\We_{\Th'\cup\Th'^\bot}\We_\Th=\We_{\Th'^\bot}\We_{\Th'}\We_\Th=\We_{\Th'^\bot}\We_\Th$.
Set $w:=w_2'\ti{w}$. We have shown that $\Th\supseteq \Th'$, and the vertices of $S(\hat{n}_{w_1'\varepsilon(R(\Th'))w_2'}u v_{(w_2)^{-1}\La})$ 
are of the form
\begin{eqnarray*}
  w_1' w \La \quad\mb{with}\quad w\in\We_{\Th'^\bot}\We_\Th  \quad\mb{such that}\quad \mb{}^\Th(w^{-1}w_2')\leq w_2\;\;.
\end{eqnarray*} 
$w_1\La$ is the biggest resp. smallest vertex of $(u_\eps)^* v_{w_1\La}$ for $\eps=+$ resp. $\eps=-$. Using 
(b) we conclude that there exists an element $w\in\We_{\Th'^\bot}\We_\Th $ such that 
$\mb{}^\Th(w^{-1}w_2')\leq w_2$ and
\begin{eqnarray*}
    w_1\La \;\left\{\begin{array}{cccc}
    \succeq & w_1' w \La &\mb{for} &\eps= +\\
    \preceq & w_1' w \La &\mb{for} & \eps=- 
            \end{array}\right.\;\;.
\end{eqnarray*}
These inequalities are equivalent to 
\begin{eqnarray*}
    w_1\;=\;w_1^\Th\; \left\{\begin{array}{cccc}
    \leq & (w_1' w)^\Th &\mb{for} &\eps= +\\
    \geq & (w_1' w)^\Th &\mb{for} & \eps=- 
            \end{array}\right.\;\;.
\end{eqnarray*} 
Obviously there also exists an element $w\in \We_{\Th'^\bot}$, which satisfies the required inequalities..\vspace*{1ex}\\
{\bf (d)} To prepare the proof of the direction from (ii) to (i) of 1a) and 2), we need the following formula: Let $w_1,w_1'\in\We$ such that $w_1\leq w_1'$ for $\eps=+$, 
and $w_1\geq w_1'$ for $\eps=-$. Let $w_2\in \mb{}^{\Th\cup\Th^\bot}\We$. By using equation (\ref{clBBG}) for the relative closures of the Bruhat and Birkhoff cells of $G$, and formula ($\delta$) stated in the part 'The minimal 
and formal Kac-Moody group $G$ and $G_f$, the monoids $\GD$ and $\GfD$' of the section 'Preliminaries' we find
\begin{eqnarray*}
  B^\eps w_1\ve{R(\Th)}w_2 B \;\,\subseteq\;\,  B^\eps w_1 B\ve{R(\Th)}w_2 B \;\,\subseteq \;\,\overline{B^\eps w_1' B}\ve{R(\Th)}w_2 B \\ 
  \;\,\subseteq \;\,\overline{B^\eps w_1' B\ve{R(\Th)}w_2 B }\;\,=\;\, 
  \overline{B^\eps w_1' \ve{R(\Th)}w_2 \underbrace{w_2^{-1}U_{\Th^\bot}w_2}_{\subseteq U^+} B }\;\,=\;\,
  \overline{B^\eps w_1' \ve{R(\Th)}w_2 B }\;\;.
\end{eqnarray*}
{\bf (e)} Now we show the direction from (ii) to (i) of 1a) and 2): By using (d) and formulas (\ref{wRwinv}), (\ref{WeDF2}) 
we find 
\begin{eqnarray*}
   (w_2')^{-1}w \ve{R(\Th)}T w_2\ &=& ((w_2')^{-1}w)^\Th \ve{R(\Th)}T w_2\;\,\subseteq\;\, \overline{B w_2^{-1}\ve{R(\Th)}w_2 B}\\
   &=& \overline{B \ve{w_2^{-1}R(\Th)} B}\;\;. 
\end{eqnarray*}
The closure $\overline{B}$ is a monoid. It contains $B$. It contains the closure $\overline{T}=\TD$, in 
particular it contains the elements $e(R)$, $R\in\RkX$. Therefore the last set of the preceeding formula is 
contained in $\overline{B}$. From this follows
\begin{eqnarray*}
   w_1'\ve{R(\Th')}w \ve{R(\Th)}w_2 T\;\,=\;\, w_1'\ve{R(\Th')}w_2'(w_2')^{-1}w \ve{R(\Th)}w_2 T\\ \;\,\in\;\,w_1'\ve{R(\Th')}w_2'\overline{B}
   \;\,\subseteq\;\,  \overline{B^\eps w_1'\ve{R(\Th')}w_2' B}\;\;.
\end{eqnarray*}
On the other hand we get by using formula (\ref{WeDF1}):
\begin{eqnarray*}
 w_1' \ve{R(\Th')} w \ve{R(\Th)} w_2 T \;\,=\;\, w_1' w \ve{R(\Th')} \ve{R(\Th)} w_2 T \;\,=\;\, w_1' w  \ve{R(\Th)} w_2 T\;\;.  
\end{eqnarray*}
Therefore 
\begin{eqnarray}\label{inc1}
  B^\eps w_1' w  \ve{R(\Th)} w_2 B\;\,\subseteq \;\,\overline{B^\eps w_1'\ve{R(\Th')}w_2' B}\;\;.
\end{eqnarray}
Due to the first inequalities of 1a) (ii) and 2) (ii) we get by using (d) and formula (\ref{WeDF2}) once more:
\begin{eqnarray}\label{inc2}
   B^\eps w_1 \ve{R(\Th)} w_2 B\;\,\subseteq \;\,\overline{B^\eps (w_1' w)^\Th \ve{R(\Th)}w_2 B}
   \;\,=\;\,\overline{B^\eps w_1' w \ve{R(\Th)}w_2 B}\;\;.
\end{eqnarray}
From the inclusions (\ref{inc1}) and (\ref{inc2}) follows (i).\\
\End
By using the explicit description of $\leq_{++}$, $\leq_{--}$, and $\leq_{-+}$ given in the last theorem, we now can 
show that these relations are order relations.
\begin{Theorem}\label{BO4}
The relations $\leq_{++}$, $\leq_{--}$, and $\leq_{-+}$ are order relations on $\WeD$, which extend the Bruhat order on $\We$.
\end{Theorem}
\Proof 
It remains to show that these relations are antisymmetric.
Let $\hat{w}= w_1 e(R(\Th)) w_2 $, $\hat{w}' =w'_1 e(R(\Th')) w'_2\in\WeD$ be elements of $\WeD$ in normal form I, and let 
$\hat{w}\leq_{++}\hat{w}'$, $\hat{w}'\leq_{++}\hat{w}$. Due to the last theorem we have $\Th=\Th'$, and there exist elements 
$w,\ti{w}\in\We_{\Th^\bot}$ such that 
\begin{eqnarray}
   w_1\leq (w_1' w)^\Th  &,&  w_1'\leq (w_1\ti{w})^\Th \;, \label{asym1.1} \\
   w_2\geq \mb{}^\Th(w^{-1}w_2')  &,& w_2'\geq \mb{}^\Th(\ti{w}^{-1}w_2)\;.\label{asym2.1}
\end{eqnarray}
Remark (2) following Proposition \ref{BO2} implies 
$\We^\Th\We_{\Th^\bot}=\We^{\Th\cup\Th^\bot}\We_{\Th^\bot}=\We^\Th$.
Therefore the inequalities of (\ref{asym1.1}) and (\ref{asym2.1}) are equivalent to
\begin{eqnarray}
   w_1\leq w_1' w  &,& w_1'\leq w_1\ti{w}\;, \label{asym1.2}\\
   w_2\geq w^{-1}w_2'  &,& w_2'\geq \ti{w}^{-1}w_2 \;. \label{asym2.2}
\end{eqnarray}
Since $w_2,w_2'\in \mb{}^{\Th\cup\Th^\bot}\We$ and $w^{-1},\ti{w}^{-1}\in\We_{\Th^\bot}\subseteq \We_{\Th\cup\Th^\bot}$ from 
(\ref{asym2.2}) follows 
\begin{eqnarray*}
  l(w_2)\geq l(w^{-1}w_2')=l(w^{-1})+l(w_2')\geq l(w_2')\geq l(\ti{w}^{-1}w_2)=l(\ti{w}^{-1})+l(w_2)\geq l(w_2)\;. 
\end{eqnarray*}
This implies $l(w_2)=l(w_2')$, and inserting in this chain of inequalities we find $l(w^{-1})=l(\ti{w}^{-1})=0$. Therefore $w^{-1}=\ti{w}^{-1}=1$. 
Inserting in (\ref{asym1.2}), (\ref{asym2.2}), we get 
\begin{eqnarray}
   w_1\leq w_1'  &,& w_1'\leq w_1 \;,\label{asym1.3}\\
   w_2\geq w_2'  &,& w_2'\geq w_2\;.\label{asym2.3}
\end{eqnarray} 
Due to the antisymmetry of the Bruhat order of $\We$ we get $w_1=w_1'$ and $w_2=w_2'$. \vspace*{1ex}\\
The antisymmetry of $\leq_{-+}$ can be proved similarly, only '$\leq$' has to be exchanged by '$\geq$' in 
(\ref{asym1.1}), (\ref{asym1.2}), (\ref{asym1.3}).\vspace*{1ex}\\
Because the inverse map $\mb{}^{inv}:\WeD\to\WeD$ is an isomorphism of the relations $\leq_{++}$ and $\leq_{--}$, 
also the relation $\leq_{--}$ is antisymmetric.\\
\End
For $\La\in P^+$ and $\la,\mu\in\We\La$ fix elements $v_\la\in L(\La)_\la\setminus\{0\}$, 
$v_\mu\in L(\La)_\mu\setminus\{0\}$ and set 
\begin{eqnarray}
 g_{\la\mu}&:=& f_{v_\la v_\mu}
\end{eqnarray}
for short. By using the antisymmetry of the extended Bruhat orders just proved we show:
\begin{Theorem}\label{BO6} Let $\hat{w}\in\WeD$. The Bruhat and Birkhoff cells $B\hat{w} B$, $B^-\hat{w} B^-$, 
and $B^-\hat{w} B$ are principal open in their closures, i.e., for   
\begin{eqnarray*}
   (\eps,\delta)= \left\{\begin{array}{l}
   (+,+) \\
   (-,-) \\
   (-,+) 
\end{array}\right\} \; \mb{ let }\; \hat{w}=w_1 \ve{R(\Th)}w_2   \;\mb{ be }\; \left\{\begin{array}{l}
   \mb{normal form I }, \\
   \mb{normal form II }, \\
   \mb{normal form I or II }, 
\end{array}\right.  
\end{eqnarray*}
and let $\La \in F_\Th\cap P^+$. Then
\begin{eqnarray}\label{sho}
       \overline{B^\eps\hat{w} B^\delta } \;\cap\; D_{\widehat{G}}(g_{w_1 \La \, w_2^{-1}\La })&=& B^\eps \hat{w} B^\delta \;\;.
\end{eqnarray}
\end{Theorem}
\Proof 
{\bf 1)} Let $\hat{w} =w_1 \ve{e(R(\Th ))}w_2 $ be in normal form I. Let $\eps\in\{+,-\}$ and $\delta=+$.\\
We first show the inclusion '$\subseteq$' of (\ref{sho}). 
Because of the Bruhat and Birkhoff decompositions of $\GD$, any element of $\GD$ can be written in the form
$u_\eps \hat{n}_{\hat{w}'}u$ with $u_\eps\in U^\eps$, $u\in U$, and $\hat{n}_{\hat{w}'}\in\ND$ belonging to 
$\hat{w}'\in\WeD$.\\ 
If $u_\eps \hat{n}_{\hat{w}'}u\in \overline{B^\eps\hat{w}B}$, then we also have $B^\eps\hat{w}'B\subseteq \overline{B^\eps\hat{w}B}$. By definition 
$\hat{w}'\leq_{++}\hat{w}$ for $\eps=+$, and $\hat{w}'\geq_{-+}\hat{w}$ for $\eps=-$.\\ 
In part (c) of the proof of Theorem \ref{BO3} we have seen that 
\begin{eqnarray*}
   g_{w_1\La w_2^{-1}\La}(u_\eps \hat{n}_{\hat{w}'}u)\neq 0
\end{eqnarray*} 
implies $\Th\supseteq \Th'$ and there exists an element $w\in\We_{\Th'^\bot}$ such that 
$\mb{}^\Th(w^{-1}w_2')\leq w_2$ and
\begin{eqnarray*}
    w_1 \left\{\begin{array}{cccc}
    \leq & (w_1' w)^\Th &\mb{for} &\eps= +\\
    \geq & (w_1' w)^\Th &\mb{for} & \eps=- 
            \end{array}\right.\;\;.
\end{eqnarray*} 
Due to Theorem \ref{BO3} itself, from this follows $\hat{w}'\geq_{++}\hat{w}$ for $\eps=+$, and 
$\hat{w}'\leq_{-+}\hat{w}$ for $\eps=-$.\\
Due to the last theorem the relations $\leq_{++}$, $\leq_{-+}$ are antisymmetric. Therefore we get 
$\hat{w}'=\hat{w}$.\vspace*{1ex}\\
To show the inclusion '$\supseteq$' of (\ref{sho}) let $u_\eps\in U^\eps$, $u\in U$, and let $\hat{n}_{\hat{w}}\in\ND$ belong to 
$\hat{w}\in\WeD=\ND/T$. Write $\hat{n}_{\hat{w}}$ in the form $\hat{n}_{\hat{w}}=n_{w_1}e(R(\Th))n_{w_2}$ with $n_{w_1}, n_{w_2}\in N$ 
belonging to $w_1, w_2\in \We$. We have
\begin{eqnarray}\label{gho}
  g_{w_1\La\, w_2^{-1}\La}(u_\eps \hat{n}_{\hat{w}}u )\;\,=\;\,   
  \kB{n_{w_1}^* u_\eps^*  v_{w_1\La}} { e(R(\Th))n_{w_2}u v_{w_2^{-1}\La} }
\end{eqnarray}
The vertices of $S(n_{w_2}u v_{w_2^{-1}\La})$ are of the form $w_2\ti{w}\La$ where $\ti{w}\in\We^{\Th}$ and $\ti{w}\leq w_2^{-1}$. Furthermore $\La$ is a vertex.\\
Now we want to determine which verticees of this form are also verticees of 
$S(e(R(\Th))n_{w_2}u v_{w_2^{-1}\La})$. Because of $\La$ is contained in the relative interior of $R(\Th)$ we 
find 
\begin{eqnarray*}
  w_2\ti{w}\La\in R(\Th) \iff w_2\ti{w} R(\Th)\subseteq R(\Th)  \iff w_2\ti{w} R(\Th) = R(\Th) \\ \iff 
  w_2\ti{w}\in \We_{\Th\cup\Th^\bot}\iff \ti{w}\in w_2^{-1}\We_{\Th\cup\Th^\bot}
\end{eqnarray*}
In this case we get $w_2^{-1}\leq \ti{w}$ because $w_2^{-1}$ is a minimal coset representative of $\We_{\Th\cup\Th^\bot}$. Because of the antisymmetry of the Bruhat order of $\We$ this implies $w_2^{-1}= \ti{w}$. Therefore 
$\La$ 
is the only vertex of $S(e(R(\Th))n_{w_2}u v_{w_2^{-1}\La})$ and we have
\begin{eqnarray*} 
 e(R(\Th))n_{w_2}u v_{w_2^{-1}\La}\;\,=\;\, n_{w_2} v_{w_2^{-1}\La}\;\;.
\end{eqnarray*} 
Inserting in equation (\ref{gho}) we find
 \begin{eqnarray*}
  g_{w_1\La\, w_2^{-1}\La}(u_\eps \hat{n}_{\hat{w}}u ) &=&  \kB{n_{w_1}^* u_\eps^*  v_{w_1\La}} { n_{w_2} v_{w_2^{-1}\La} } \\  
                                                     &=& \kB{n_{w_1}^*   v_{w_1\La}} { n_{w_2} v_{w_2^{-1}\La} } \;\,\neq\;\; 0
\end{eqnarray*}
{\bf 2)} By applying the Chevalley involution $*:\GD\to\GD$ on (\ref{sho}) for $\eps\in\{ +,-\}$, and $\delta=+$ we 
find the equation
\begin{eqnarray*}
   \overline{B^-(\hat{w})^{inv} B^{-\eps}} \;\cap \;D_{\widehat{G}}(g^*)\;\,=\;\,B^- (\hat{w})^{inv} B^{-\eps}\;\;,
\end{eqnarray*}
from which follow the remaining statements of the theorem by Remark (1) following Proposition \ref{BO2}.\\ 
\End
The next two theorems give product decompositions of $B\hat{w}B$, $B^-\hat{w}B^-$, and $B^-\hat{w}B$, 
$\hat{w}\in \WeD$, as principal open sets in their closures. 
These generalize the product decompositions of sets $BwB=U_w\cdot wT\cdot U$, and $B^-wB=(U^w)^- \cdot wT\cdot U $, 
$w\in\We$, given in \cite{KP1} for Kac-Moody groups. The decomposition of  $B\hat{w}B$, $\hat{w}\in\WeD$, is 
nearly similar to the corresponding decomposition of a Bruhat cell of the wonderful compactification of a 
semisimple algebraic group in \cite{Re3}.
\begin{Theorem}\label{BO7}\mb{}\\
{\bf 1a)} Let $\hat{w}=w_1\ve{R(\Th)}w_2$ be an element of $\WeD$ in normal form I, and let $\hat{n}$ be a corresponding element of $\ND$.
Then
\begin{eqnarray*}
    B\,\hat{w}\,B &=& U_{w_1}\, \hat{n}\, (w_2^{-1}T^\Th w_2)\,(U\cap w_2^{-1}U^\Th w_2)\;\;.
\end{eqnarray*}
{\bf b)} Let $\hat{w}=w_1\ve{R(\Th)}w_2$ be an element of $\WeD$ in normal form II, and let $\hat{n}$ be a corresponding element of $\ND$.
Then
\begin{eqnarray*}
   B^-\,\hat{w}\,B^- &=& (U^- \cap w_1 (U^\Th)^-  w_1^{-1}) \,(w_1T^\Th w_1^{-1}) \, \hat{n}\, U_{w_2^{-1}}^-\;\;.
\end{eqnarray*}
{\bf 2)}  Let $\hat{w}=w_1\ve{R(\Th)}w_2$ be an element of $\WeD$ in normal form I or in normal form II, and let $\hat{n}$ be a corresponding 
element of $\ND$. Then 
\begin{eqnarray*}
    B^-\,\hat{w}\,B &=& (U^-\cap w_1 (U^\Th)^-w_1^{-1})\,\hat{n} \,(w_2^{-1}T^\Th w_2)\,(U\cap w_2^{-1}U^\Th w_2)\\
                    &=& (U^-\cap w_1 (U^\Th)^-w_1^{-1})\,(w_1 T^\Th w_1^{-1})\, \hat{n} \,(U\cap w_2^{-1}U^\Th w_2)\;\;.
\end{eqnarray*}
\end{Theorem}
{\bf Remark:} It is easy to see that 
\begin{eqnarray*}
  U_{w_1} &=& U\cap w_1 U^- w_1^{-1}\;\,=\;\, U\cap w_1 (U^\Th)^- w_1^{-1}\;\;,\\
  U_{w_2^{-1}} &=& U\cap w_2^{-1} U^- w_2\;\,=\;\, U\cap w_2^{-1} (U^\Th)^- w_2\;\;.
 \end{eqnarray*}
In this sense the formulas of 1) are symmetric in the first and last factor.\vspace*{1ex}\\
\Proof 
We only have to show the statements of the theorem which involve normal form I. Then the statements which involve normal form II follow by applying 
the Chevalley involution $*:\GD\to\GD$, together with the Remark (1) following Proposition \ref{BO2}.\vspace*{1ex}\\
Write $\hat{n}$ in the form $\hat{n}=n_1e(R(\Th))n_2$ with $n_1, n_2\in N$ corresponding to $w_1\in\We^\Th$, $w_2\in\mb{}^{\Th\cup\Th^\bot}\We$. For the following transformations we make use of the formulas ($\beta$), 
($\gamma$), and ($\delta$) stated in the part ``The minimal and formal Kac-Moody group $G$ and $G_f$, the 
monoids $\GD$ and $\GfD$'' of the section ``Preliminaries''.\vspace*{1ex}\\
{\bf 1)} By using the first equation of ($\delta$) we find
\begin{eqnarray*}
    B\hat{w}B &=& U_{w_1} U^{w_1} n_1 e(R(\Th)) T n_2 U \;\,=\;\, 
    U_{w_1} n_1  \underbrace{w_1^{-1} U^{w_1}w_1}_{\subseteq U}  e(R(\Th)) T n_2 U  \\
    &\subseteq & U_{w_1} n_1  e(R(\Th)) U_{\Th^\bot} T n_2 U 
    \;\,=\;\, U_{w_1} n_1 e(R(\Th)) T n_2 \underbrace{w_2^{-1}U_{\Th^\bot}w_2 }_{\subseteq U} U \\
    &=& U_{w_1} n_1 e(R(\Th)) T n_2 U \;\,=\;\, U_{w_1} n_1 e(R(\Th)) T n_2 U_{w_2^{-1}}U^{w_2^{-1}} \\
    &=& U_{w_1} n_1 e(R(\Th)) w_2 U_{w_2^{-1}}w_2^{-1} T n_2 U^{w_2^{-1}}\;\;.
\end{eqnarray*}
$w_2 U_{w_2^{-1}}w_2^{-1}$ is generated by the root groups
\begin{eqnarray*}
  U_\beta & \mb{where}& \beta \in \nrW\;,\; w_2^{-1}\beta\in\prW\;\;.
\end{eqnarray*}
Since $w_2^{-1}\in\We^{\Th\cup\Th^\bot}$ these root groups coincide with the root groups 
\begin{eqnarray*}
  U_\beta & \mb{where}& \beta \in \nrW\setminus \We_{\Th\cup\Th^\bot}\Mklz{\al_i}{i\in (\Th\cup\Th^\bot)} \;,\; w_2^{-1}\beta\in\prW\;\;,
\end{eqnarray*}
which are contained in $(U^{\Th\cup\Th^ \bot})^-$. By using formula (\ref{wRwinv}) we find
\begin{eqnarray*}
    B\hat{w}B \;\,\subseteq\;\, U_{w_1} n_1 e(R(\Th)) T n_2 U^{w_2^{-1}}\;\;,
\end{eqnarray*}
and the reverse inclusion is obvious.
We have
\begin{eqnarray*}
  U^{w_2^{-1}} &\subseteq & w_2^{-1}U w_2 \;\,=\,\;w_2^{-1}U_\Th w_2 \ltimes w_2^{-1}U^\Th w_2\;\;.
\end{eqnarray*}
Because of $w_2^{-1}\in\We^{\Th\cup\Th^\bot}$, the group $w_2^{-1}U_\Th w_2 $ is contained in $U$. Clearly it is also contained in 
$w_2^{-1}U w_2 $. Therfore it is contained in $U^{w_2^{-1}}$.\\
An element $x\in w_2^{-1}U w_2 $ can be written in the form $x=p_1(x)p_2(x)$ with $p_1(x)\in w_2^{-1}U_\Th w_2 $ and $p_2(x)\in w_2^{-1}U^\Th w_2$.
Obviously $p_2(U^{w_2^{-1}})\supseteq U^{w_2^{-1}}\cap w_2^{-1}U^\Th w_2$. The reverse inclusion follows because for $x\in U^{w_2^{-1}}$ we have
\begin{eqnarray*}
   w_2^{-1}U^\Th w_2\;\ni \;p_2(x)\;=\;p_1(x)^{-1}x \;\in\; U^{w_2^{-1}}\;\;.
\end{eqnarray*}
We get
\begin{eqnarray*}
   p_2(U^{w_2^{-1}}) \;\,=\;\,U^{w_2^{-1}}\cap w_2^{-1} U^\Th w_2 \;\,=\;\, U\cap w_2^{-1} U w_2 \cap w_2^{-1} U^\Th w_2 
   \;\,=\;\, U\cap w_2^{-1} U^\Th w_2\;\;.
\end{eqnarray*}
By using this equation and two times formula ($\beta$) we find
\begin{eqnarray*}
    B\hat{w}B &=& U_{w_1} n_1 e(R(\Th)) T n_2 U^{w_2^{-1}}
    \;\;=\;\, U_{w_1} n_1 e(R(\Th))\underbrace{w_2 p_1(U^{w_2^{-1}})w_2^{-1}}_{\subseteq U_\Th} T n_2 p_2(U^{w_2^{-1}}) \\
      &=& U_{w_1} n_1 e(R(\Th))T n_2 (U\cap w_2^{-1} U^\Th w_2) 
    \;\,=\;\, U_{w_1} n_1 e(R(\Th))T^\Th n_2 (U\cap w_2^{-1} U^\Th w_2) \\
     &=& U_{w_1} n_1 e(R(\Th)) n_2 (w_2^{-1}T^\Th w_2)  (U\cap w_2^{-1} U^\Th w_2) \;\;.
\end{eqnarray*}
{\bf 2)} By using the first equation of ($\delta$) we get
\begin{eqnarray*}
    B^-\hat{w}B &=& (U^{w_1})^- U_{w_1}^-  n_1 e(R(\Th)) T n_2 U \;\,=\;\, 
           (U^{w_1})^- n_1 \underbrace{w_1^{-1}U_{w_1}^-  w_1}_{\subseteq U^+} e(R(\Th)) T n_2 U \\
             &\subseteq & (U^{w_1})^- n_1  e(R(\Th)) U_{\Th^\bot}T n_2 U 
             \;\,=\;\, (U^{w_1})^- n_1  e(R(\Th))T n_2 \underbrace{w_2^{-1}U_{\Th^\bot}w_2}_{\subseteq U} U \\
             &=& (U^{w_1})^- n_1  e(R(\Th))T n_2 U \;\;.
\end{eqnarray*}
The reverse inclusion is obvious. In the same way as before we get
\begin{eqnarray*}
    B^-\hat{w}B  &=& (U^{w_1})^- n_1  e(R(\Th))T n_2 (U\cap w_2^{-1} U^\Th w_2 )\\
    &=& (U^{w_1})^- n_1  e(R(\Th)) n_2 (w_2^{-1}T^\Th w_2) (U\cap w_2^{-1} U^\Th w_2 )\;\;. 
\end{eqnarray*}
Treating the first factor in a similar way we find
\begin{eqnarray*}
    B^-\hat{w}B &=& (U^-\cap w_1(U^\Th)^- w_1^{-1}) n_1  e(R(\Th)) n_2 (w_2^{-1}T^\Th w_2) (U\cap w_2^{-1} U^\Th w_2 )\\
    &=& (U^-\cap w_1(U^\Th)^- w_1^{-1}) (w_1 T^\Th w_1^{-1}) n_1  e(R(\Th)) n_2  (U\cap w_2^{-1} U^\Th w_2 )\;\;. 
\end{eqnarray*}
\End
We equip $B^\eps\hat{w}B^\delta$ with its coordinate ring $\FK{B^\eps\hat{w}B^\delta}$ as a principal open set in 
its closure, $(\eps,\delta)\in\{(++), \,(--),\,(-+)\}$, $\hat{w}\in\WeD$.\\
Recall that the torus  $T$ of the Kac-Moody group can be described by the following isomorphism of 
groups:\vspace*{1.5ex}\\
 \hspace*{8em} $\begin{array}{ccc}
 H\otimes_\Z\F^\times &\to & \;\;\:T \\
 \sum_{i=1}^{2n-l}h_i\otimes s_i\;\, &\mapsto &\, \prod_{i=1}^{2n-l} t_{h_i}(s_i) 
\end{array}\;$.\vspace*{1.5ex}\\
The group algebra $\FK{P}$ of the lattice $P$ can be identified with the classical coordinate ring on $T$,
identifying $\sum c_\la e_\la\in \FK{P}$ with the function on $T$ defined by
\begin{eqnarray*}
  \left(\sum_\la c_\la e_\la\right)\left(\prod_{i=1}^{2n-l}t_{h_i}(s_i)\right) \;\,:=\,\;\sum_\la \,c_\la \,
  \prod_{i=1}^{2n-l}(s_i)^{\la(h_i)} &\;\,,\;\,&   (s_i\in\F^\times )\;\;.
\end{eqnarray*}
Similarly, for $J\subseteq I$ and $w\in\We$, the classical coordinate ring of the torus $w T^J w^{-1}$, 
where 
\begin{eqnarray*}
    T^J:=\{\,\prod_{i=1,\,\ldots,\, 2n-l \atop i\notin J}t_{h_i}(s_i)\,\mid\, s_i\in\F^\times \,\}\;\;,
\end{eqnarray*} 
is given by the group algebra $\FK{w P^J}$, where 
\begin{eqnarray*}
   P^J:= \Z\mb{-span}\Mklz{\La_i}{i=1,\,\ldots,\,2n-l,\; i\notin J}
\end{eqnarray*}
(In general the classical coordinate rings of these tori do not coincide with the restriction of the coordinate ring 
$\FK{\GD}$ onto these tori. It is possible to show that $w T^J w^{-1}$ is principal open in its closure, and the 
classical coordinate ring of $w T^J w^{-1}$ is the coordinate ring of this principal open set. But we do not need 
this for the following considerations.) 
\begin{Theorem}\label{BO8}\mb{}\\
{\bf 1a)} Let $\hat{w}=w_1\ve{R(\Th)}w_2$ be an element of $\WeD$ in normal form I, and let $\hat{n}$ be a corresponding element of $\ND$. The map 
\begin{eqnarray*}
   m:\, U_{w_1} \,\times\, w_2^{-1}T^\Th w_2 \,\times\,  (U\cap w_2^{-1}U^\Th w_2)   \;\to\;  B\hat{w}B \;\;,
\end{eqnarray*}
defined by $m(u,t,\ti{u}):= u \hat{n}t\ti{u}$, is an isomorphism. \vspace*{1ex}\\
{\bf b)} Let $\hat{w}=w_1\ve{R(\Th)}w_2$ be an element of $\WeD$ in normal form II, and let $\hat{n}$ be a corresponding element of $\ND$. The map 
\begin{eqnarray*}
   m:\, (U^- \cap w_1 (U^\Th)^-  w_1^{-1}) \,\times\, w_1 T^\Th w_1^{-1} \,\times\,  U_{w_2^{-1}}^-   \;\to\;  B^-\hat{w} B^- \;\;,
\end{eqnarray*}
defined by $m(u,t,\ti{u}):= u t \hat{n} \ti{u}$, is an isomorphism.\vspace*{1ex}\\
{\bf 2)}  Let $\hat{w}=w_1\ve{R(\Th)}w_2$ be an element of $\WeD$ in normal form I or in normal form II, and let $\hat{n}$ be a corresponding 
element of $\ND$. The map 
\begin{eqnarray*}
   m:\, (U^- \cap w_1 (U^\Th)^-  w_1^{-1}) \,\times\, w_2^{-1} T^\Th w_2 \,\times\,  (U\cap w_2^{-1}U^\Th w_2)   \;\to\;  B^-\hat{w} B \;\;,
\end{eqnarray*}
defined by $m(u,t,\ti{u}):= u \hat{n}t\ti{u}$, is an isomorphism. 
\end{Theorem}
\Proof We only have to show the statements of the theorem which involve normal form I. Then the statements which involve normal form II follow 
easily by using the Chevalley involution $*:\GD\to\GD$, together with Remark (3) following Definition \ref{BO1} 
and Remark (1) following Proposition \ref{BO2}.\vspace*{1ex}\\ 
Due to the last theorem the maps $m$ of 1a) and 2) are surjective. We show that the corresponding comorphisms $m^*$ are 
well defined and surjective. This is sufficient, because the surjectivity of the maps $m$ imply the injectivity of the 
comorphisms $m^*$, and the surjectivity of the comorphisms $m^*$ imply the injectivity of the maps $m$.\vspace*{1ex}\\ 
{\bf To 1a)}: To show that the comorphism 
\begin{eqnarray*}
   m^*:  \FK{B\hat{w}B} \;\to\; \FK{U_{w_1}} \,\otimes\, \FK{w_2^{-1}P^\Th} \,\otimes\,  \FK{U\cap w_2^{-1}U^\Th w_2}    
\end{eqnarray*}
is well defined, write $\hat{n}$ in the form $\hat{n}=n_1e(R(\Th))n_2$ with $n_1, n_2\in N$ corresponding to $w_1\in\We^\Th$, $w_2\in\mb{}^{\Th\cup\Th^\bot}\We$.\\
Fix an element $\La\in F_\Th$, and fix $v_1\in L(\La)_{w_1\La}\setminus\{0\}$, $v_2\in L(\La)_{w_2^{-1}\La}\setminus\{0\}$. 
For $u\in U_{w_1}$, $t\in w_2^{-1}T^\Th w_2$, and $\ti{u}\in (U\cap w_2^{-1}U^\Th w_2)$ we have
\begin{eqnarray*}
  f_{v_1v_2}(m(u,t,\ti{u})) \;\,=\;\, \kB{v_1}{u n_1e(R(\Th))n_2 t\ti{u}v_2} \;\,=\;\,\kB{u^* v_1}{n_1e(R(\Th))n_2 t\ti{u}v_2}\;\;. 
\end{eqnarray*}
Because of $\ti{u}\in w_2^{-1}U^\Th w_2$ we get $\ti{u}v_2=v_2$. Because of 
\begin{eqnarray*}
  u^*\;\,\in\;\,U_{w_1}^*\;\,=\;\,(U\cap w_1 U^- w_1^{-1})^*\;\,=\;\, U^-\cap w_1 U w_1^{-1}
\end{eqnarray*}
we get $u v_1=v_1$. Therefore 
\begin{eqnarray}\nonumber
  f_{v_1v_2}(m(u,t,\ti{u})) &=& \kB{v_1}{n_1 e(R(\Th)) n_2 t v_2} \;\,=\;\, \kB{v_1}{n_1 e(R(\Th)) n_2 v_2} e_{w_2^{-1}\La}(t)\\
                            &=& \underbrace{\kB{v_1}{n_1 n_2 v_2}}_{\neq 0} e_{w_2^{-1}\La}(t)\;\;.\label{fv1v2}
\end{eqnarray}
Now let $N\in P^+$ and $v,w\in L(N)$. Choose $\kBl$-dual bases of $L(N)$ by choosing $\kBl$-dual bases 
\begin{eqnarray*}
  (a_{\tau i})_{i=1,\,\ldots,\, m_\tau} &,&  (b_{\tau i})_{i=1,\,\ldots,\, m_\tau}
\end{eqnarray*} 
of $L(N)_\tau$ for every $\tau\in P(N)$. For $u\in U_{w_1}$, $t\in w_2^{-1}T^\Th w_2$, and $\ti{u}\in (U\cap w_2^{-1}U^\Th w_2)$ we have
\begin{eqnarray*}
  f_{vw}(m(u,t,\ti{u})) &=& \kB{v}{un_1e(R(\Th))n_2t\ti{u}w} \\
   &=& \sum_{\tau,\,i} \kB{v}{un_1e(R(\Th))n_2 t a_{\tau i}}\kB{b_{\tau i}}{\ti{u}w}\\
   &=& \sum_{\tau,\,i\atop w_2\tau\in R(\Th)} \kB{v}{un_1n_2 a_{\tau i}}e_\tau(t)\kB{b_{\tau i}}{\ti{u}w}\;\;.
\end{eqnarray*}
A summand of this sum is nonzero at the most if $w\neq 0$, and if $\tau$ is bigger than an element of $supp(w)$, which is only possible for 
finitely many non-zero summands. Furthermore $R(\Th)\cap P = X\cap P^\Th$. Taking into account (\ref{fv1v2}) we therefore get for $p\in\Nn$:
\begin{eqnarray}
    \lefteqn{m^*\left(\frac{f_{vw}}{(f_{v_1 v_2})^p}\res{B\hat{w}B}\right)}  \label{wd} \\ 
   &=&  \sum_{\tau,\,i\atop \tau\in w_2^{-1}P^\Th} \frac{1}{\kB{v_1}{n_1n_2v_2}^p} \,f_{v\,n_1n_2a_{\tau i}}\res{U_{w_1}}\,\otimes\,
         e_{\tau-p(w_2^{-1})\La}\,\otimes\, f_{b_{\tau i}w}\res{U\cap w_2^{-1}U^\Th w_2}\;\;. \nonumber
\end{eqnarray}
Because $\FK{\GD}$ is spanned by the matrix coefficients $f_{vw}$, $v,w\in L(N)$, $N\in P^+$, we have shown that 
the comorphism $m^*$ is well defined.\vspace*{1ex}\\
To show the surjectivity of $m^*$, it is sufficient to find elements of $\FK{B\hat{w}B}$, which are mapped onto a system of generators of 
$\FK{U_{w_1}} \,\otimes\, \FK{w_2^{-1}P^\Th} \,\otimes\,  \FK{U\cap w_2^{-1}U^\Th w_2}$.\vspace*{1ex}\\
a) Let $v\in L(\La)$. For $u\in U_{w_1}$, $t\in w_2^{-1}T^\Th w_2$, and $\ti{u}\in (U\cap w_2^{-1}U^\Th w_2)$ we have
\begin{eqnarray*}
  \kB{v_1}{u n_1 e(R(\Th))n_2 t \ti{u} v}  \;\,=\;\, \kB{u^* v_1}{n_1 e(R(\Th))n_2 t \ti{u} v} \;\,=\;\, \kB{v_1}{n_1 e(R(\Th))n_2 t \ti{u} v}\\
          \;\,=\;\, \kB{t n_2^* e(R(\Th))n_1^*v_1}{\ti{u} v}\;\,=\;\,  \kB{ t n_2^* n_1^* v_1}{\ti{u} v}
           \;\,=\;\, \kB{n_2^* n_1^* v_1}{\ti{u} v} e_{w_2^{-1}\La}(t)\;\;.
\end{eqnarray*}
Taking into accout (\ref{fv1v2}) we therefore get 
\begin{eqnarray}\label{sa)}
  m^*\left(\frac{f_{v_1 v}}{f_{v_1 v_2}}\res{B\hat{w}B}\right)
  &=& \underbrace{\frac{1}{\kB{v_1}{n_1n_2v_2}}}_{\neq 0} \, 1\,\otimes\,1\,\otimes\, f_{n_2^* n_1^* v_1 \,v}\res{U\cap w_2^{-1}U^\Th w_2} \;\;.
\end{eqnarray} 
We have $n_1^* v_1\in L(\La)_\La\setminus\{0\}$. Due to Theorem 5.6 of \cite{M1} the coordinate ring $\FK{U^\Th}$ is 
generated by the functions $f_{n_1^* v_1\,v}\res{U^\Th}$, $v\in L(\La)$. Therefore its restriction 
$\FK{w_2 U w_2^{-1}\cap U^\Th}$ is generated by 
\begin{eqnarray*}
  f_{n_1^* v_1\, n_2 v}\res{w_2 U w_2^{-1}\cap U^\Th} &\quad,\quad& v\in L(\La)\;\;.
\end{eqnarray*} 
It is easy to check that we get an isomorphism 
\begin{eqnarray*}
\phi: U\cap w_2^{-1}U^\Th w_2 \to w_2 U w_2^{-1}\cap U^\Th 
\end{eqnarray*}
by $\phi(u) :=n_2 u n_2^{-1}$. For $u\in U\cap w_2^{-1}U^\Th w_2$ we have
\begin{eqnarray*}
  \phi^*(f_{n_1^* v_1\,n_2 v}\res{w_2 U w_2^{-1}\cap U^\Th})(u) &=& \kB{n_1^* v_1}{n_2 u n_2^{-1} n_2 v} \;\,= \;\, 
                                                                     \kB{ n_2^* n_1^* v_1}{u v}\;\;.
\end{eqnarray*}
Therefore the functions 
\begin{eqnarray*}
\phi^*(f_{n_1^* v_1\,n_2 v}\res{w_2 U w_2^{-1}\cap U^\Th})\;\,= \;\,f_{n_2^* n_1^* v_1\,v}\res{U\cap w_2^{-1}U^\Th w_2 }\qquad,\qquad v\in L(\La)\;,
\end{eqnarray*}
which appear in (\ref{sa)}), generate the coordinate ring $\FK{U\cap w_2^{-1}U^\Th w_2}$.
\vspace*{1ex}\\
b) Let $N\in\overline{F_\Th}\cap P$ and $\ti{v}_1\in L(N)_{w_1 N}\setminus\{0\}$, $\ti{v}_2\in L(N)_{w_2^{-1} N}\setminus\{0\}$. Let $p\in\Nn$. Similar to (\ref{fv1v2}) we find
\begin{eqnarray}\label{sb)}
  m^*\left(\frac{f_{\ti{v}_1 \ti{v}_2}}{(f_{v_1 v_2})^p}\res{B\hat{w}B}\right)
     &=& \underbrace{\frac{\kB{\ti{v}_1}{n_1n_2 \ti{v}_2}}{\kB{v_1}{n_1n_2v_2}^p}}_{\neq 0} \, 1\,\otimes\,e_{w_2^{-1}(N - p\La)}\,\otimes\,1\;\;.
\end{eqnarray}
It is easy to check that $(\overline{F_\Th}\cap P)-\Nn\La=P^\Th$. Therefore the functions $e_{w_2^{-1}(N - m\La)}$, $N\in\overline{F_\Th}\cap P$, 
$m\in\Nn$, span the coordinate ring $\FK{w_2^{-1}P^\Th}$.\vspace*{1ex}\\
c) Let $v\in L(\La)$. For $u\in U_{w_1}$, $t\in w_2^{-1}T^\Th w_2$, and $\ti{u}\in (U\cap w_2^{-1}U^\Th w_2)$ we have
\begin{eqnarray*}
  \kB{v}{u n_1 e(R(\Th))n_2 t \ti{u} v_2} &=& \kB{v}{u n_1 e(R(\Th))n_2 t v_2} \;\,=\;\, \kB{v}{u n_1 e(R(\Th))n_2 v_2} e_{w_2^{-1}\La}(t)\\
                                          &=& \kB{v}{u n_1 n_2 v_2} e_{w_2^{-1}\La}(t)\;\;.
\end{eqnarray*}
Taking into account (\ref{fv1v2}) we therefore get
\begin{eqnarray}\label{sc)}
  m^*\left(\frac{f_{v v_2}}{f_{v_1 v_2}}\res{B\hat{w}B}\right)
        &=& \underbrace{\frac{1}{\kB{v_1}{n_1 n_2 v_2}}}_{\neq 0} \, f_{v\,n_1 n_2 v_2} \res{U_{w_1}}\,\otimes\,1\,\otimes\,1\;\;.
\end{eqnarray}
$U_{w_1}=U\cap w_1 U^- w_1^{-1}$ is generated by the root groups $U_\al$, $\al\in \prW\cap w_1 \nrW $. Since 
$w_1\in\We^\Th$ we have $\prW\cap w_1 \nrW = \prW\cap w_1 (\nrW\setminus \We_\Th\Th)$. This implies
\begin{eqnarray*}
   U_{w_1} \;\,=\;\, U\cap w_1 U^- w_1^{-1}\;\,\subseteq \;\,U\cap w_1 (U^\Th)^- w_1^{-1}\;\;,
\end{eqnarray*}
and the reverse inclusion is obvious. We have $n_2 v_2\in L(\La)_\La\setminus\{0\}$. Due to Theorem 5.6 of 
\cite{M1} the coordinate ring $\FK{(U^\Th)^-}$ is generated by the functions $f_{v \,n_2 v_2}\res{(U^\Th)^-}$, 
$v\in L(\La)$. Therefore its restriction $\FK{w_1^{-1}U w_1 \cap (U^\Th)^-}$ is generated by 
\begin{eqnarray*}
    f_{n_1^*v \,n_2 v_2}\res{w_1^{-1}U w_1 \cap (U^\Th)^-} &\quad,\quad & v\in L(\La)\;\;.
\end{eqnarray*}
It is easy to check that we get an isomorphism 
\begin{eqnarray*}
\phi: U \cap w_1 (U^\Th)^- w_1^{-1}\to w_1^{-1}U w_1 \cap (U^\Th)^-
\end{eqnarray*}
by $\phi(u) := n_1^{-1}u n_1$. For $u\in U_{w_1}= U \cap w_1 (U^\Th)^- w_1^{-1}$ we have
\begin{eqnarray*}
      \phi^*(f_{n_1^* v \,n_2 v_2}\res{w_1^{-1}U w_1 \cap (U^\Th)^-})(u)\;\,=\;\, 
      \kB{n_1^* v}{n_1^{-1}u n_1 n_2 v_2}\;\,=\;\, \kB{v}{u n_1 n_2 v_2}
\end{eqnarray*}
Therefore the functions 
\begin{eqnarray*}
  \phi^*(f_{n_1^* v \,n_2 v_2}\res{w_1^{-1}U w_1 \cap (U^\Th)^-})\;\,=\;\,
f_{v\, n_1 n_2 v_2}\res{U_{w_1}}\qquad, \qquad v\in L(\La)\;,
\end{eqnarray*}
which appear in (\ref{sc)}) generate the coordinate ring $\FK{U_{w_1}}$.\vspace*{1ex}\\
{\bf To 2):} We use the same notations as in the first part of this proof.
For $u\in U^-\cap w_1 (U^\Th)^- w_1^{-1}$ we have $u^*\in w_1 U^\Th w_1^{-1}$. Therefore also $u^*v_1=v_1$.\\ 
Completely parallel to the first part of this proof we get formula (\ref{wd}) with $B\hat{w}B$ replaced by $B^-\hat{w}B$, and $U_{w_1}$ replaced by $U^-\cap w_1 (U^\Th)^- w_1^{-1}$.
In particular it follows that the comorphism 
\begin{eqnarray*}
   m^*:  \FK{B^-\hat{w}B} \;\to\; \FK{U^- \cap w_1 (U^\Th)^-  w_1^{-1}} \,\otimes\, \FK{w_2^{-1}P^\Th} \,\otimes\,  \FK{U\cap w_2^{-1}U^\Th w_2}    
\end{eqnarray*}
is well defined.\vspace*{1ex}\\
Also completely parallel to a) and b) of the first part of this proof we get the formulas (\ref{sa)}) and 
(\ref{sb)}) with $B\hat{w}B$ replaced by $B^-\hat{w}B$. Therefore we have found elements of 
$\FK{B^-\hat{w}B}$, which are mapped onto a system of generators of 
$1\,\otimes\, 1 \,\otimes\,  \FK{U\cap w_2^{-1}U^\Th w_2}$, and $1\,\otimes\, \FK{w_2^{-1}P^\Th} \,\otimes\,  1 $.\vspace*{1ex}\\
Completely parallel to the proof of formula (\ref{sc)}), we get formula (\ref{sc)}) with
$B\hat{w}B$ replaced by $B^-\hat{w}B$, and $U_{w_1}$ replaced by $U^-\cap w_1 (U^\Th)^- w_1^{-1}$.
An easy modification of the corresponding argument of a) of the first part of the proof shows that the 
functions $f_{v\, n_1 n_2 v_2}\res{U^-\cap w_1 (U^\Th)^- w_1^{-1}}$, 
$v\in L(\La)$, generate the coordinate ring $\FK{U^-\cap w_1 (U^\Th)^- w_1^{-1}}$.\\
\End
Next we want to show that the orbits are irreducible. For this we first state a theorem of \cite{M2}.
Equip $\GfD$ with a coordinate ring by identifying with the coordinate ring of $1\tr \GfD \subseteq \Spm\FK{\GD}$ 
via the injective map $\GfD\to1\tr\GfD$, $x\mapsto 1\tr x$. Denote the Zariski closure of $M\subseteq \GfD$ by 
$\overline{\overline{M}}$. To use later note that due to Theorem 16 of \cite{M2} this closure coincides with the 
closure denoted by $\overline{\overline{M}}$ in \cite{M2}. Now Theorem 4 of \cite{M2} gives:
\begin{Theorem}\label{BO9} Let $D_1$, $D_2$ be irreducible subgroups of $G_f$, and 
$x\in\GD$. The $D_1\times D_2$-orbit of the element $1 \tr x\in \Spm\FK{G}$ is irreducible.
\end{Theorem}
There is a similar theorem for $\GD$:
\begin{Theorem}\label{BO10}
Let $D_1$, $D_2$ be irreducible subgroups of $G$, and let $x\in\GD$. The $D_1\times D_2$-orbit $D_1 ^* x D_2$ is 
irreducible.
\end{Theorem}
The proof of this theorem can be extracted from a part of the proof of Theorem \ref{BO9}. For the convenience of the reader we sketch the proof.\\ 
\Proof Let $x\in\GD$ and $Or:=D_1^* xD_2$ its $D_1\times D_2$-orbit. Let $A_1$ and $A_2$ be closed subsets of $\GD$, such that 
$Or\subseteq A_1\cup A_2$. We have to show $Or\subseteq A_1$ or $Or\subseteq A_2$.\vspace*{1ex}\\
Let $d_1\in D_1$. The map $\gamma_{d_1}:D_2\to \GD$ defined by $\gamma_{d_1}(d_2):= d_1^* x d_2$, 
$d_2\in D_2$, is a morphism.\\ 
Similarly, for $d_2\in D_2$, the map $\delta_{d_2}:D_1\to \GD$ defined by $\delta_{d_2}(d_1) :=  d_1^* x d_2$, $d_1\in D_1$, is a 
morphism.\vspace*{1ex}\\
Let $d_1\in D_1$. Because of $\gamma_{d_1}(D_2)\subseteq Or$, we have $\gamma_{d_1}^{-1}(A_1)\cup \gamma_{d_1}^{-1}(A_2)=D_2$. Furthermore 
$\gamma_{d_1}^{-1}(A_1)$ and $\gamma_{d_1}^{-1}(A_2)$ are closed. Because of the 
irreducibility of $D_2$ we get $\gamma_{d_1}^{-1}(A_1)=D_2$ or $\gamma_{d_1}^{-1}(A_2)=D_2$.\vspace*{1ex}\\  
Therefore the sets
\begin{eqnarray*}
  B_1 &:=& \Mklz{d_1\in D_1}{\gamma_{d_1}^{-1}(A_1)=D_2}\;\;,\\
  B_2 &:=& \Mklz{d_1\in D_1}{\gamma_{d_1}^{-1}(A_2)=D_2}
\end{eqnarray*}
satisfy $B_1\cup B_2 =D_1$.\\
Note that for $d_1\in D_1$ and $d_2\in D_2$ we have $\gamma_{d_1}(d_2)=\delta_{d_2}(d_1)$. The set $B_1$ is closed, because of 
\begin{eqnarray*}
  B_1 &=& \Mklz{d_1\in D_1}{ \gamma_{d_1}(d_2)\in A_1 \;\mb{ for all }\; d_2\in D_2}\;\;=\;\;\bigcap_{d_2\in D_2}\underbrace{\delta_{d_2}^{-1}(A_1)}_{closed}\;\;.
\end{eqnarray*}
Similarly, the set $B_2$ is closed.
Because of the irreducibility of $D_1$ we get $B_1=D_1$ or $B_2=D_1$, which is equivalent to $Or\subseteq A_1$ or $Or\subseteq A_2$.\\
\End
\begin{Prop}\label{BO11} $B$ and $B^-$ are irreducible.
\end{Prop}
\Proof
The Chevalley involution $*:\GD\to\GD$ is an isomorphism, which maps $B$ onto $B^-$. Therefore it is sufficient to show that 
$B$ is irreducible.\vspace*{1ex}\\  
Due to Corollary \ref{BO6} we have $D_{\overline{B}}(g_{\La\La}\res{\overline{B}})=B$. Due to 
Theorem \ref{BO8} the multiplication map 
\begin{eqnarray*}
  m: U\times T\to D_{\overline{B}}(g_{\La\La}\res{\overline{B}})
\end{eqnarray*}
is an isomorphism. Here $U$ is irreducible because its coordinate ring is a symmetric algebra, 
compare \cite{KP2}, Lemma 4.3. The torus $T$, which is equipped with its classical coordinate ring, is 
irreducible. Therefore also the principal open set $D_{\overline{B}}(g_{\La\La}\res{\overline{B}})$ is 
irreducible, which implies that $\overline{B}$ is irreducible, which implies that $B$ is irreducible.\\
\End
Combining the last proposition with Theorem \ref{BO10} we get:
\begin{Cor}\label{BO12} The Bruhat and Birkhoff cells $B\hat{w}B$, $B^-\hat{w}B^-$, and $B^-\hat{w}B$ are 
irreducible for every $\hat{w}\in\WeD$.
\end{Cor}
At last we consider the Birkhoff cells of $\Spm\FK{G}$. They have similar properties as the Birkhoff cells 
of $\GD$: 
\begin{Theorem}\label{BO13} Let $\hat{w}=w_1\ve{R(\Th)}w_2$ be an element of $\WeD$ in normal form I or in normal form II, and let $\hat{n}$ be a corresponding element of $\ND$.\\
The closure of the Birkhoff cell $B_f\tr \hat{w} B_f$ is given by 
\begin{eqnarray}\label{Spmcl}
  \overline{B_f\tr \hat{w}B_f}^{Spm} &=& \bigcup_{\hat{w}'\in\widehat{\cal W} \atop \hat{w}'\geq_{-+} \hat{w}} B_f\tr \hat{w}' B_f\;¸\;.
\end{eqnarray}
The Birkhoff cell $B_f\tr \hat{w} B_f$ is irreducible. It is principal open in its closure, i.e., 
\begin{eqnarray}\label{Spmpo}
       \overline{B_f\tr \hat{w} B_f} \;\cap\; D_{Specm \,\F\,[\widehat{G}]}(g_{w_1 \La \, w_2^{-1}\La })&=& B_f\tr \hat{w} B_f \;\;.
\end{eqnarray} 
Equip $B_f\tr \hat{w} B_f$ with its coordinate ring as a principal open set in its closure. Equip the torus 
$w_2^{-1} T^\Th w_2$ with its classical coordinate ring. Then the map
\begin{eqnarray}\label{Spmiso}
   m:\, (U_f \cap w_1 U^\Th_f  w_1^{-1}) \,\times\, w_2^{-1} T^\Th w_2 \,\times\,  (U_f\cap w_2^{-1}U_f^\Th w_2)   \;\to\;  B_f\tr\hat{w} B_f \;\;,
\end{eqnarray}
defined by $m(u,t,\ti{u}):=u\tr\hat{n}t\ti{u}$, is an isomorphism.
\end{Theorem}
\Proof
We first show the equation (\ref{Spmcl}): It is easy to adapt step (c), $\eps=-$, of the proof of 
Theorem \ref{BO3} to show the inclusion '$\subseteq$'. To show the reverse inclusion, note that for $M\subseteq \GD$ 
we have 
\begin{eqnarray*}
   1\,\tr\overline{M}\;\subseteq\; \overline{1\,\tr M}^{Spm}\;\;.
\end{eqnarray*}
Let $\hat{w}'\in\WeD$ such that $\hat{w}'\geq_{-+} \hat{w}$, then
\begin{eqnarray*}
  B\tr\hat{w}'B &=& 1\,\tr B^-\hat{w}'B\;\,\subseteq\,\; 1\,\tr \overline{B^-\hat{w} B}\;\,\subseteq\;\, \overline{1\,\tr B^-\hat{w} B}^{Spm} 
  \,\;=\,\; \overline{B\,\tr \hat{w}'B}^{Spm}\\
   &\subseteq & \overline{B_f\,\tr \hat{w} B_f}^{Spm}\;\;.
\end{eqnarray*}
Since $B_f\times B_f$ acts on $\Spm\FK{G}$ by morphisms, we find by applying $B_f\times B_f$ to this inclusion:
\begin{eqnarray*}
   B_f\tr\hat{w}'B_f \;\subseteq\; \overline{B_f\,\tr \hat{w} B_f}^{Spm}\;\;.
\end{eqnarray*}
It is easy to adapt the proof of Theorem \ref{BO6} to show equation (\ref{Spmpo}).
\vspace*{1ex}\\
Because of Theorem \ref{BO9}, to show the irreducibility of $B_f\hat{w} B_f$, it is sufficient to show the 
irreducibility of $B_f\subseteq \widehat{G_f}$.\\
The relative topology on B induced by $\widehat{G_f}$ is the same as the relative topology on $B$ 
induced by $\GD$. Due to Proposition \ref{BO11} $B$ is irreducible. Therefore it is sufficient to show
\begin{eqnarray}\label{BfB}
         \overline{\overline{B_f}}  &=&  \overline{\overline{B}}\;\;.
\end{eqnarray}
Due to \cite{M2}, Theorem 9 (1) we have $\overline{\overline{U}}=U_f$. Due to Proposition 1 of \cite{M2} 
$\overline{\overline{B}}$ is a monoid. Because $T$ and $U$ are contained in $B$ we find 
$B_f= T U_f\subseteq \overline{\overline{B}}$, from which follows the inclusion '$\subseteq$' in (\ref{BfB}). 
The reverse inclusion is obvious.\vspace*{1ex}\\
Due to Theorem 9 (2) of \cite{M2} we have $\overline{\overline{U^\Th}}= U_f^\Th$. Therefore the coordinate 
ring $\FK{U_f^\Th}$ is isomorphic to $\FK{U^\Th}$ by the restriction map. Due to part (a) of the proof of 
Theorem 21 in \cite{M2} we have $e(R(\Th))(U_f)_\Th=e(R(\Th))$. Using these results it is not difficult to 
adapt the proof of Theorem \ref{BO7} and Theorem \ref{BO8} to show that the map $m$ in (\ref{Spmiso}) is an 
isomorphism.\\
\End
%
%
%
\section{Extensions of the length function\label{SL}} 
%
%
%
We first define three extensions of the length function $l:\We\to\Nn$ of the Weyl group to functions 
$l_{++},l_{--}:\WeD\to\Z$ and $l_{-+}:\WeD\to\Nn$ of the Weyl monoid by using the normal forms of the 
elements of $\WeD$. The length function $l_{++}$ is similar to the 
length function of a Renner monoid given in \cite{Re2} up to additive constants on the orbits given by the 
action of the product of the Weyl group on the Renner monoid. The functions $l_{++}$, $l_{--}$ are allowed to 
take positive and negative values. It is not possible to make $l_{++}$, $l_{--}$ positive on the 
$\We\times\We$-orbits of $\WeD$ by adding constants, because in general these functions can take arbitrary 
positive and negative values on such orbits. 
\begin{Def}\label{L1}\mb{}\\
{\bf 1a)} Let $\hat{w}=w_1\ve{R(\Th)}w_2$ be an element of $\WeD$ in normal form I. Set
\begin{eqnarray*}
    l_{++}(\hat{w}) &:=& l(w_1)-l(w_2)\;\in\;\Z\;\;.
\end{eqnarray*}
{\bf b)} Let $\hat{w}=w_1\ve{R(\Th)}w_2$ be an element of $\WeD$ in normal form II. Set
\begin{eqnarray*}
    l_{--}(\hat{w}) &:=& -l(w_1)+l(w_2)\;\in\;\Z\;\;.
\end{eqnarray*}
{\bf 2)}  Let $\hat{w}=w_1\ve{R(\Th)}w_2$ be an element of $\WeD$ in normal form I or in normal form II. Set
\begin{eqnarray*}
    l_{-+}(\hat{w}) &:=& l(w_1)+l(w_2)\;\in\;\Nn\;\;.
\end{eqnarray*}
\end{Def}
{\bf Remarks:}\\ 
{\bf (1)} We have $l_{++}(\hat{w})=l_{--}(\hat{w}^{inv})$, and $l_{-+}(\hat{w})=l_{-+}(\hat{w}^{inv})$.\\
{\bf (2)} If $\hat{w}=w_1 w_2\ve{R(\Th)} w_3=w_1\ve{R(\Th)} w_2 w_3$ is an element of $\WeD$ in normal form III, then 
\begin{eqnarray*}
   l_{++}(\hat{w}) &=& \;\;\;l(w_1) +l(w_2)-l(w_3)\;\;,\\ 
   l_{--}(\hat{w}) &=&    -l(w_1) +l(w_2)+l(w_3)\;\;,\\ 
   l_{-+}(\hat{w}) &=& \;\;\;l(w_1) +l(w_2)+l(w_3)\;\;.\vspace*{0.5ex}
\end{eqnarray*}
The length function $l:\We\to\Nn$ is compatible with the Bruhat order on $\We$ and the natural order on $\Nn$. 
Similar things hold for the the extensions of the length function restricted to the $\We\times\We$-orbits of $\WeD$:
\begin{Prop} Let $\Th$ be special, and $\hat{w},\hat{w}'\in \We\ve{R(\Th)}\We$. Equip $\We\ve{R(\Th)}\We$ with the restriction of the extended 
Bruhat order $\leq_{\eps\delta}$, $(\eps\delta)\in\{(++),(--),(-+)\}$.\vspace*{1ex}\\
a)  Then $\hat{w}\lneqq_{\eps\delta}\hat{w}'$ implies $l_{\eps\delta}(\hat{w})\lneqq l_{\eps\delta}(\hat{w}')$.\vspace*{1ex}\\
b) The length of every chain joining $\hat{w}$ and $\hat{w}'$ is finite, and does not extend 
$l_{\eps\delta}(\hat{w}')-l_{\eps\delta}(\hat{w})$. In particular there exists a maximal chain joining 
$\hat{w}$ and $\hat{w}'$.
\end{Prop}
\Proof
b) follows immediately from a). We only have to prove the cases $(++),(-+)$ of a). Then the case $(--)$ follows from the case $(++)$ by using the inverse map $\mb{}^{inv}\WeD\to\WeD$, combining Remark (3) following 
Definition \ref{BO1} and Remark (1) following Definition \ref{L1}.\vspace*{1ex}\\ 
Recall that for elements $u,v\in\We$ we have
\begin{eqnarray}\label{Wu1}
 |\,l(u)-l(v)| \;\,\leq\;\, l(uv) \;\,\leq\;\, l(u)+l(v)\;\;,
\end{eqnarray}
and for $J\subseteq I$, $u_J\in\We_J$, $\mb{}^J u\in\mb{}^J\We$ we have
\begin{eqnarray}\label{Wu2}
  l(u_J\mb{}^J u)\;\,=\;\,l(u_J)+l(\mb{}^J u)\;\;.
\end{eqnarray}
Let $\hat{w}=w_1\ve{R(\Th)}w_2$, $\hat{w}'=w_1'\ve{R(\Th)}w_2'$ be in normal form I, i.e., $w_1,w_1'\in\We^\Th$ and 
$w_2,w_2'\in\mb{}^{\Th\cup\Th^\bot}\We$.\\
{\bf a)} If $\hat{w}\lneqq_{++}\hat{w}'$ then due to Theorem \ref{BO3} there exist an element $w\in\We_{\Th^\bot}$, such that
\begin{eqnarray*}
    w_1\leq (w_1'w)^\Th = w_1' w \quad\mb{ and }\quad  w_2 \geq \mb{}^\Th (w^{-1}w_2')= w^{-1}w_2' \;\;.
\end{eqnarray*} 
Furthermore we have $w_1 \neq w_1' w$ or $w_2 \neq w^{-1}w_2' $ because otherwise we get by using formula 
(\ref{WeDF1}):
\begin{eqnarray*}
  \hat{w}\;\,=\;\, w_1\ve{R(\Th)}w_2 \;\,=\;\, w_1'w\ve{R(\Th)}w^{-1}w_2'\;\,=\;\, w_1'\ve{R(\Th)}w_2' \;\,=\;\, \hat{w}'\;\;.  
\end{eqnarray*}
Now $l:(\We,\leq)\to(\Nn,\leq)$ is an order morphism. By using equation (\ref{Wu2}) and inequality (\ref{Wu1}) 
we get
\begin{eqnarray*}
 l_{++}(\hat{w}) &=&       l(w_1)-l(w_2) \;\,\lneqq\;\, l(w_1'w)-l(w^{-1}w_2') \;\,=\;\, l(w_1'w)-l(w^{-1})-l(w_2')\\
                 &\leq &   l(w_1'w w^{-1})-l(w_2')\;\,=\;\,l_{++}(\hat{w}')\;\;.
\end{eqnarray*}
{\bf b)} If $\hat{w}\lneqq_{-+}\hat{w}'$ then due to Theorem \ref{BO3} there exist an element $w\in\We_{\Th^\bot}$, 
such that
\begin{eqnarray*}
    w_1'\geq (w_1 w)^\Th = w_1 w \quad\mb{ and }\quad  w_2' \geq \mb{}^\Th (w^{-1}w_2)= w^{-1}w_2 \;\;.
\end{eqnarray*} 
Furthermore we have $w_1' \neq  w_1 w$ or $w_2' \neq  w^{-1}w_2$ because otherwise 
\begin{eqnarray*}
  \hat{w}'\;\,=\;\, w_1'\ve{R(\Th)}w_2' \;\,=\;\, w_1 w\ve{R(\Th)}w^{-1}w_2 \;\,=\;\, w_1\ve{R(\Th)}w_2 \;\,=\;\, \hat{w} \;\;.  
\end{eqnarray*}
By using that $l:(\We,\leq)\to(\Nn,\leq)$ is an order morphism, by using equation (\ref{Wu2}) and 
inequality (\ref{Wu1}) we get
\begin{eqnarray*}
 l_{-+}(\hat{w}') &=&       l(w_1')+ l(w_2') \;\,\gneqq\;\, l(w_1 w)+ l(w^{-1}w_2) \;\,=\;\, l(w_1 w)+l(w^{-1}) +l(w_2')\\
                 &\geq &   l(w_1 w w^{-1})+l(w_2)\;\,=\;\,l_{-+}(\hat{w})\;\;.
\end{eqnarray*}
\End
For a reductive algebraic monoid the length of any maximal chain of the Bruhat-Chevalley order, which joins two 
elements in an orbit of the action of the product of the Weyl group on the Renner monoid, is given by the 
difference of the length functions of these elements. This follows from the definition of the length function 
in \cite{Re2}, the algebraic description of the Bruhat Chevalley order in \cite{PePuRe}, and the Theorem in 
Section 7 of \cite{Re1}, which states that if an element covers another then the dimension difference of the 
corresponding Bruhat cells is one. The algebraic geometric proof of this theorem can not be generalized to the 
Kac-Moody setting.\\
A combinatorial proof working for a Renner monoid and the Weyl monoid is not easy to find, although some 
particular cases are not difficult to  treat. In a separate article this will be investigated further.
\vspace*{1ex}\\   
For an element $w\in\We$ and $i\in I$ the length of $\sigma_i w$ and $w$, and the length of $w\sigma_i$ and 
$w$ are related in the following way:
\begin{eqnarray}
  l(\sigma_i w) &=& \left\{\begin{array}{lcl} 
                                l(w)+1  & \mb{for} &  w^{-1}\al_i\in\prW\\ 
                                l(w)-1  & \mb{for} &  w^{-1}\al_i\in\nrW
                          \end{array}\right.\;\;, \label{sw}\\
  l(w\sigma_i) &=& \left\{\begin{array}{lcl} 
                                l(w)+1  & \mb{for} &  w\al_i\in\prW\\ 
                                l(w)-1  & \mb{for} &  w\al_i\in\nrW
                          \end{array}\right. \;\;. \label{ws}
\end{eqnarray}
The next theorem gives the generalization for the extended length functions. To cut short our notation we set 
$\We_J J:=\We_J \Mklz{\al_i}{i\in J}$, $J\subseteq I$.
\begin{Theorem} {\bf 1)} Let $\hat{w}=w_1\ve{R(\Th)}w_2$ be an element of $\WeD$ in normal form I. Let $i\in I$. We have:
\begin{eqnarray*}
    l_{++}(\sigma_i \hat{w}) &=& \left\{\begin{array}{lcl}
                              l_{++}(\hat{w})+1 & \mb{ for } & w_1^{-1}\al_i\in \prW\setminus \We_{\Th}\Th  \\
                              l_{++}(\hat{w})   & \mb{ for } & w_1^{-1}\al_i\in \We_{\Th}\Th \\
                              l_{++}(\hat{w})-1 & \mb{ for } & w_1^{-1}\al_i\in \nrW\setminus \We_{\Th}\Th\\
                            \end{array}\right. \\
    l_{--}(\sigma_i \hat{w}) &=& \left\{\begin{array}{lcl}
                              l_{--}(\hat{w})+1 & \mb{ for } & w_1^{-1}\al_i\in \left(\nrW\setminus \We_{\Th \cup \Th^\bot}(\Th\cup\Th^\bot)\right) 
\cup \left(\prW\cap \We_{\Th^\bot}\Th^\bot\right)\\
                              l_{--}(\hat{w})   & \mb{ for } & w_1^{-1}\al_i\in \We_{\Th}\Th \\
                              l_{--}(\hat{w})-1 & \mb{ for } & w_1^{-1}\al_i\in \left(\prW\setminus \We_{\Th \cup \Th^\bot}(\Th\cup\Th^\bot)\right) 
    \cup \left(\nrW\cap \We_{\Th^\bot}\Th^\bot\right)
                            \end{array}\right.\\
 l_{-+}(\sigma_i \hat{w}) &=& \left\{\begin{array}{lcl}
                              l_{-+}(\hat{w})+1 & \mb{ for } & w_1^{-1}\al_i\in \prW\setminus \We_{\Th}\Th  \\
                              l_{-+}(\hat{w})   & \mb{ for } & w_1^{-1}\al_i\in \We_{\Th}\Th \\
                              l_{-+}(\hat{w})-1 & \mb{ for } & w_1^{-1}\al_i\in \nrW\setminus \We_{\Th}\Th\\
                            \end{array}\right.
\end{eqnarray*}
In all three cases the extended lengths of $\sigma_i\hat{w}$ and $\hat{w}$ are equal if and only if the 
elements $\sigma_i\hat{w}$ and $\hat{w}$ are equal.\vspace*{1ex}\\ 
{\bf 2)} Let $\hat{w}=w_1\ve{R(\Th)}w_2$ be an element of $\WeD$ in normal form II. Let $i\in I$. We have:
\begin{eqnarray*}
    l_{++}(\hat{w}\sigma_i) &=& \left\{\begin{array}{lcl}
                              l_{++}(\hat{w})+1 & \mb{ for } & w_2\al_i\in \left(\nrW\setminus \We_{\Th \cup \Th^\bot}(\Th\cup\Th^\bot)\right) 
\cup \left(\prW\cap \We_{\Th^\bot}\Th^\bot\right)\\
                              l_{++}(\hat{w})   & \mb{ for } & w_2\al_i\in \We_{\Th}\Th \\
                              l_{++}(\hat{w})-1 & \mb{ for } & w_2\al_i\in \left(\prW\setminus \We_{\Th \cup \Th^\bot}(\Th\cup\Th^\bot)\right) 
    \cup \left(\nrW\cap \We_{\Th^\bot}\Th^\bot\right)
                            \end{array}\right. \\
  l_{--}(\hat{w}\sigma_i) &=& \left\{\begin{array}{lcl}
                              l_{++}(\hat{w})+1 & \mb{ for } & w_2\al_i\in \prW\setminus \We_{\Th}\Th  \\
                              l_{++}(\hat{w})   & \mb{ for } & w_2\al_i\in \We_{\Th}\Th \\
                              l_{++}(\hat{w})-1 & \mb{ for } & w_2\al_i\in \nrW\setminus \We_{\Th}\Th\\
                            \end{array}\right.\\
l_{-+}(\hat{w}\sigma_i) &=& \left\{\begin{array}{lcl}
                              l_{-+}(\hat{w})+1 & \mb{ for } & w_2\al_i\in \prW\setminus \We_{\Th}\Th  \\
                              l_{-+}(\hat{w})   & \mb{ for } & w_2\al_i\in \We_{\Th}\Th \\
                              l_{-+}(\hat{w})-1 & \mb{ for } & w_2\al_i\in \nrW\setminus \We_{\Th}\Th\\
                            \end{array}\right.  
\end{eqnarray*}
In all three cases the extended lengths of $\hat{w}\sigma_i$ and $\hat{w}$ are equal if and only if the 
elements $\hat{w}\sigma_i$ and $\hat{w}$ are equal. 
\end{Theorem}
\Proof We only have to show the statements 2) of the theorem, which involve normal form II. Then the statements 
1), which involve normal form I, follow easily by using Remark (1) following Definition \ref{L1} and Remark (1) 
following Proposition \ref{BO2}. We make use of the formulas (\ref{WeDF1}) and (\ref{WeDF2}) without further mentioning.\vspace*{1ex}\\
Let $\hat{w}= w_1 \ve{R(\Th)}w_2$ be in normal form II, i.e., $w_1\in \We^{\Th\cup\Th^\bot}$, $w_2\in\mb{}^\Th\We$. Let $w_2=bc$ with 
$b\in\We_{\Th^\bot}$ and $c \in \mb{}^{\Th\cup\Th^\bot}\We$.\vspace*{1ex}\\
{\bf a)} We have:
\begin{eqnarray*}
     \hat{w}\sigma_i=\hat{w} &\iff& \ve{R(\Th)}w_2\sigma_i w_2^{-1}=\ve{R(\Th)}\\
  &\iff& \sigma_{w_2\al_i } = w_2\sigma_i w_2^{-1}\in\We_\Th   \;\,\iff\;\, w_2\al_i\in\We_\Th\Th
\end{eqnarray*}
{\bf b)} Let $w_2\al_i\in\We_{\Th^\bot}\Th^\bot$, which is equivalent to $c\al_i\in\We_{\Th^\bot}\Th^\bot$. Then 
\begin{eqnarray*}
   \hat{w}\sigma_i\;\,=\;\,w_1 b\ve{R(\Th)}c \sigma_i \;\,=\;\, w_1 b \ve{R(\Th)}\sigma_{c\al_i }c\;\;=\;\, 
   w_1 b \sigma_{c\al_i }\ve{R(\Th)}c\;\;.
\end{eqnarray*}
Note that the expression on the right is in normal form III with $b \sigma_{c\al_i }$ as middle part. Now $c\al_i$ 
can be written in the form $c\al_i=\sum_{j\in\Th}n_j\al_j$ with either all $n_j\in\Nn$ or all $n_j\in-\Nn$. 
Applying $c^{-1}$ we get
\begin{eqnarray*}
 \al_i &=& \sum_{j\in\Th} n_j c^{-1}\al_j\;\;.
\end{eqnarray*}
Because of $c^{-1}\in\We^{\Th\cup\Th^\bot}$ we have $c^{-1}\al_j\in\prW\subseteq \Nn\mb{-span}\{\al_1,\,\ldots,\,\al_n\}$ 
for all $j\in\Th$. This implies that there exists an index $j\in\Th$ such that $\al_i=c^{-1}\al_j$. In particular 
$c\al_i$ is a simple root. By using Remark (2) following Definition \ref{L1} and equation (\ref{ws}) we find
\begin{eqnarray*}
   l_{++}(\hat{w}\sigma_i) -l_{++}(\hat{w})\;\,=\;\,l_{--}(\hat{w}\sigma_i) -l_{--}(\hat{w})\;\,=\;\,l_{-+}(\hat{w}\sigma_i) -l_{-+}(\hat{w})\\
    \;\,=\;\,l(b\sigma_   {c\al_i })-l(b)
    \;\;=\;\,\left\{\begin{array}{lcl}
              1  & \mb{for} & w_2\al_i = bc\al_i\in \prW  \\
              -1  & \mb{for} & w_2\al_i= bc\al_i\in \nrW
         \end{array}\right.\;\;.
\end{eqnarray*}
{\bf c)} Let $w_2\al_i\in\rW\setminus \We_{\Th\cup\Th^\bot}(\Th\cup\Th^\bot)$, which is equivalent to 
$c\al_i\in\rW\setminus \We_{\Th\cup\Th^\bot}(\Th\cup\Th^\bot)$.\\
We first show $c\sigma_i\in\mb{}^{\Th\cup\Th^\bot}\We$: Because of $c^{-1}\in\We^{\Th\cup\Th^\bot}$ we have
\begin{eqnarray*}
             c^{-1}\al_j\in\prW \quad\mb{ for all }\quad j\in\Th\cup\Th^\bot \;\;.
\end{eqnarray*}
Due to our assumption on $c\al_i$, from this follows
\begin{eqnarray*}
  c^{-1}\al_j\in\prW\setminus\{\al_i\}\;\subseteq\; \prW\setminus\{\al_i\}\cup\{-\al_i\} \quad\mb{ for all }\quad  j\in\Th\cup\Th^\bot\;\;.  
\end{eqnarray*}
Applying $\sigma_i$ we get 
\begin{eqnarray*}
\sigma_i c^{-1}\al_j\in\prW \quad\mb{ for all }\quad   j\in\Th\cup\Th^\bot\;,
\end{eqnarray*}
which implies $\sigma_i c^{-1}\in \We^{\Th\cup\Th^\bot}$, or equivalent $c\sigma_i \in\mb{}^{\Th\cup\Th^\bot}\We$. \vspace*{1ex}\\
Now $\hat{w}\sigma_i=w_1\ve{R(\Th)}b c\sigma_i$ is in normal form III with $c\sigma_i$ as last part. By using 
Remark (2) following Definition \ref{L1} and equation (\ref{ws}) we find
\begin{eqnarray*}
l_{++}(\hat{w}\sigma_i)-l_{++}(\hat{w}) \;\,=\;\, -l_{--}(\hat{w}\sigma_i)+l_{--}(\hat{w})\;\,=\;\,-l_{-+}(\hat{w}\sigma_i)+l_{-+}(\hat{w})\\
        =\;\,-l(c\sigma_i)+l(c) \;\,=\;\,
          \left\{ \begin{array}{lcl}
                +1 & \mb{for} & c\al_i\in\nrW\setminus \We_{\Th\cup\Th^\bot}(\Th\cup\Th^\bot)\\
                -1 & \mb{for} & c\al_i\in\prW\setminus \We_{\Th\cup\Th^\bot}(\Th\cup\Th^\bot)
           \end{array}\right.\;\;.
\end{eqnarray*}
Here $c\al_i\in\rW^\pm\setminus \We_{\Th\cup\Th^\bot}(\Th\cup\Th^\bot)$ is equivalent to 
$w_2\al_i=bc\al_i\in\rW^\pm\setminus \We_{\Th\cup\Th^\bot}(\Th\cup\Th^\bot)$.\vspace*{1ex}\\
Putting together the cases a), b), and c) part 2) of the theorem follows.\\
\End
For the monoid of matrices over a finite field, and more general for a reductive algebraic monoid a 
generalization of one of the Tits axioms for BN-pairs has been given in \cite{So}, \cite{Re2} by using the length function.
The next theorem gives a generalization of some of the Tits axioms for groups with twinned BN-pairs, compare 
\cite{Ti}, for the monoid $\GD$:  
\begin{Theorem} Let $i\in I$, let $\hat{w}\in\WeD$. Let $\eps\in\{+,-\}$. We have:
\begin{eqnarray*}
  (B^\eps\sigma_i B^\eps)\,(B^\eps\hat{w}B^\eps) &=& \left\{\begin{array}{lcl}
                B^\eps\hat{w}B^\eps & \mb{ if } & l_{\eps\eps}(\sigma_i\hat{w}) =l_{\eps\eps}(\hat{w})\\
                B^\eps\sigma_i \hat{w}B^\eps & \mb{ if } & l_{\eps\eps}(\sigma_i\hat{w}) =l_{\eps\eps}(\hat{w})+1\\
                B^\eps\sigma_i \hat{w}B^\eps \;\dot{\cup}\;B^\eps\hat{w}B^\eps & \mb{ if } & l_{\eps\eps}(\sigma_i\hat{w}) =l_{\eps\eps}(\hat{w})-1\\  
                                 \end{array}\right.\;\;.
\end{eqnarray*}
\begin{eqnarray*}
  (B^\eps\hat{w}B^\eps)\,(B^\eps\sigma_i B^\eps)  &=& \left\{\begin{array}{lcl}
                   B^\eps\hat{w}B^\eps & \mb{ if } & l_{\eps\eps}(\hat{w}\sigma_i) =l_{\eps\eps}(\hat{w})\\
                   B^\eps \hat{w}\sigma_i B^\eps & \mb{ if } & l_{\eps\eps}(\hat{w}\sigma_i) =l_{\eps\eps}(\hat{w})+1\\
                   B^\eps \hat{w}\sigma_i B^\eps \;\dot{\cup}\;B^\eps\hat{w}B^\eps & \mb{ if } & l_{\eps\eps}(\hat{w}\sigma_i) =l_{\eps\eps}(\hat{w})-1\\  
                                 \end{array}\right.\;\;.
\end{eqnarray*}
We have:
\begin{eqnarray*}
  (B^{-\eps} \sigma_i B^{-\eps} )\,(B^{-\eps} \hat{w}B^\eps) &=& \left\{\begin{array}{lcl}
            B^{-\eps}\hat{w}B^\eps & \mb{ if } & l_{\eps\eps}(\sigma_i\hat{w}) =l_{\eps\eps}(\hat{w})\\
            B^{-\eps}\sigma_i \hat{w}B^\eps & \mb{ if } & l_{\eps\eps}(\sigma_i\hat{w}) =l_{\eps\eps}(\hat{w})-1\\
            B^{-\eps}\sigma_i \hat{w}B^\eps \;\dot{\cup}\;B^{-\eps}\hat{w}B^\eps & \mb{ if } & l_{\eps\eps}(\sigma_i\hat{w}) =l_{\eps\eps}(\hat{w})+1\\  
                                 \end{array}\right.\;\;.
\end{eqnarray*}
\begin{eqnarray*}
  (B^\eps \hat{w} B^{-\eps} )\,(B^{-\eps} \sigma_i B^{-\eps}) &=& \left\{\begin{array}{lcl}
               B^\eps\hat{w}B^{-\eps} & \mb{ if } & l_{\eps\eps}(\hat{w}\sigma_i) =l_{\eps\eps}(\hat{w})\\
               B^\eps\hat{w}\sigma_i B^{-\eps} & \mb{ if } & l_{\eps\eps}(\hat{w}\sigma_i) =l_{\eps\eps}(\hat{w})-1\\
               B^\eps\hat{w}\sigma_i B^{-\eps} \;\dot{\cup}\;B^\eps\hat{w}B^{-\eps} & \mb{ if } & l_{\eps\eps}(\hat{w}\sigma_i) =l_{\eps\eps}(\hat{w})+1\\  
                                 \end{array}\right.\;\;.
\end{eqnarray*}
\end{Theorem}
\Proof For the following transformations we use the formulas ($\beta$) and ($\gamma$) stated in the part ``The minimal 
and formal Kac-Moody group $G$ and $G_f$, the monoids $\GD$ and $\GfD$'' of the section ``Preliminaries'' several 
times.\\
{\bf 1)} Let $\hat{w}=w_1\ve{R(\Th)}w_2$ be in normal form I, i.e., $w_1 \in \We^\Th$, $w_2\in\mb{}^{\Th\cup\Th^\bot}\We$.\vspace*{1ex}\\
{\bf a)} We have
\begin{eqnarray*}
   \sigma_i B \hat{w} &=& \sigma_i T U^i U_i   w_1\ve{R(\Th)} w_2 \;\,=\;\, T U^i \sigma_i U_i   w_1\ve{R(\Th)} w_2 \\ 
                      &=& T U^i \sigma_i w_1 U_{w_1^{-1}\al_i} \ve{R(\Th)} w_2\\
                      &=& \left\{\begin{array}{lcl}
                           T U^i \sigma_i w_1 \ve{R(\Th)} w_2  & for & 
                                                                w_1^{-1}\al_i\in\prW\setminus \We_{\Th\cup\Th^\bot}(\Th\cup\Th^\bot) \cup\We_\Th\Th\\
                           T U^i \sigma_i w_1 \ve{R(\Th)} U_{w_1^{-1}\al_i} w_2  & for & w_1^{-1}\al_i\in\We_{\Th^\bot}\Th^\bot 
                          \end{array} \right.\\
                      &=& \left\{\begin{array}{lcl}
                           T U^i \sigma_i \hat{w} & for & w_1^{-1}\al_i\in\prW\setminus \We_{\Th\cup\Th^\bot}(\Th\cup\Th^\bot)\cup\We_\Th\Th \\
                           T U^i \sigma_i \hat{w} U_{w_2^{-1}w_1^{-1}\al_i}   & for & w_1^{-1}\al_i\in\We_{\Th^\bot}\Th^\bot
                          \end{array} \right.\;\;.
\end{eqnarray*}
Since $w_2^{-1}\in\We^{\Th\cup\Th^\bot}$ we have $w_2^{-1}w_1^{-1}\al_i\in (\rW)^\pm$ 
if $w_1^{-1}\al_i\in\We_{\Th^\bot}\Th^\bot \cap (\rW)^\pm$. We conclude:
\begin{eqnarray*}
   B \sigma_i B \hat{w} B    &=&  B\sigma_i \hat{w} B    \quad \mb{ for }\quad  w_1^{-1}\al_i\in\left(\prW\setminus \We_\Th \Th \right)\cup\We_\Th\Th\;\;.\\
   B \sigma_i B \hat{w} B^-  &=&  B\sigma_i \hat{w} B^-  \mb{ for } w_1^{-1}\al_i\in 
                                 \left(\prW\setminus \We_{\Th\cup\Th^\bot}(\Th\cup\Th^\bot)\right) \cup 
                                 (\We_{\Th^\bot}\Th^\bot \cap \nrW)  \cup  \We_\Th\Th\;\;.
\end{eqnarray*}
{\bf b)} We have
\begin{eqnarray*}
  \sigma_i B \hat{w} &=& \sigma_i T U^i U_i  w_1\ve{R(\Th)} w_2 \;\,=\;\, T U^i \sigma_i U_i   w_1\ve{R(\Th)} w_2 \\
                     &=&  T U^i (\sigma_i U_i\sigma_i^{-1}) \sigma_i w_1\ve{R(\Th)} w_2 \;\,\subseteq\;\, 
                          T U^i (U_i\cup U_i\sigma_i U_i) \sigma_i w_1\ve{R(\Th)} w_2 \\
                     &=&  B\sigma_i \hat{w} \;\,\cup\;\, B\sigma_i(\sigma_i w_1 U_{(\sigma_i w_1)^{-1}\al_i})\ve{R(\Th)}w_2 \\
                     &=&  B\sigma_i \hat{w} \;\,\cup\;\, \left\{ \begin{array}{lcl}
                          B w_1\ve{R(\Th)}w_2 & for &  w_1^{-1}\al_i \in \nrW\setminus \We_{\Th\cup\Th^\bot}(\Th\cup\Th^\bot)\\
                          B w_1\ve{R(\Th)}U_{-w_1^{-1}\al_i} w_2 & for &  w_1^{-1}\al_i \in \We_{\Th^\bot}\Th^\bot
                          \end{array}\right.\\
                     &=&  B\sigma_i \hat{w} \;\,\cup\;\, \left\{ \begin{array}{lcl}
                          B \hat{w} & for &  w_1^{-1}\al_i \in \nrW\setminus \We_{\Th\cup\Th^\bot}(\Th\cup\Th^\bot)\\
                          B \hat{w}U_{-w_2^{-1}w_1^{-1}\al_i} & for &  w_1^{-1}\al_i \in \We_{\Th^\bot}\Th^\bot
                          \end{array}\right.\;\;.\\                                             
\end{eqnarray*}
Here $-w_2^{-1}w_1^{-1}\al_i\in (\rW)^\mp$ if $w_1^{-1}\al_i\in\We_{\Th^\bot}\Th^\bot \cap (\rW)^\pm$. Note also 
that $\sigma_i U_i\sigma_i^{-1}$ contains elements of $U_i$ and $U_i\sigma_i U_i$. We conclude:
\begin{eqnarray*}
   B \sigma_i B \hat{w} B    &=&  B\sigma_i \hat{w} B \cup  B\hat{w}B    \quad\mb{ for }\quad w_1^{-1}\al_i\in\nrW\setminus \We_\Th \Th \;\;.\\
   B \sigma_i B \hat{w} B^-  &=&  B\sigma_i \hat{w} B^- \cup  B\hat{w}B^- \quad\mb{ for }\quad w_1^{-1}\al_i\in 
                                 \left(\nrW\setminus \We_{\Th\cup\Th^\bot}(\Th\cup\Th^\bot)\right) 
                                  \cup (\We_{\Th^\bot}\Th^\bot \cap \prW) \;\;.
\end{eqnarray*}
Due to Part 1) of the last theorem, we get from a) and b): 
\begin{eqnarray*}
   B \sigma_i B \hat{w} B &=& \left\{\begin{array}{lcl}
                                           B\hat{w}B                            &\mb{if}&   l_{++}(\sigma_i\hat{w})=l_{++}(\hat{w}) \\
                                           B\sigma_i\hat{w}B                    &\mb{if}&   l_{++}(\sigma_i\hat{w})=l_{++}(\hat{w})+1 \\
                                           B\sigma_i\hat{w}B \;\dot{\cup}\; B\hat{w}B     &\mb{if}&   l_{++}(\sigma_i\hat{w})=l_{++}(\hat{w})-1 
                                          \end{array}\right.\;\;.\\
   B \sigma_i B \hat{w} B^- &=& \left\{\begin{array}{lcl}
                                           B\hat{w}B^-                              &\mb{if}&   l_{--}(\sigma_i\hat{w})=l_{--}(\hat{w}) \\
                                           B\sigma_i\hat{w}B^-                      &\mb{if}&   l_{--}(\sigma_i\hat{w})=l_{--}(\hat{w})-1 \\
                                           B\sigma_i\hat{w}B^- \;\dot{\cup}\; B\hat{w}B^-     &\mb{if}&   l_{--}(\sigma_i\hat{w})=l_{--}(\hat{w})+1 
                                          \end{array}\right.\;\;.
\end{eqnarray*}
{\bf 2)}  Let $\hat{w}=w_1\ve{R(\Th)}w_2$ be in normal form II, i.e., $w_1 \in \We^{\Th\cup\Th^\bot}$ and $w_2\in \mb{}^\Th \We$.\vspace*{1ex}\\
{\bf a)} We have
\begin{eqnarray*}
  \hat{w}B\sigma_i &=& w_1\ve{R(\Th)}w_2 U_i U^i T\sigma_i\;\,=\;\, w_1\ve{R(\Th)}w_2 U_i \sigma_i U^i T \\
  &=& w_1\ve{R(\Th)} U_{w_2\al_i}w_2\sigma_i U^i T\\
  &=& \left\{\begin{array}{lcl}
           w_1\ve{R(\Th)}w_2 \sigma_i U^i T &\mb{ for } &  w_2\al_i\in \nrW\setminus\We_{\Th\cup\Th^\bot}(\Th\cup\Th^\bot)\cup \We_\Th\Th\\
           w_1 U_{w_2\al_i}\ve{R(\Th)}w_2 \sigma_i U^i T &\mb{ for }& w_2\al_i\in\We_{\Th^\bot}\Th^\bot
      \end{array}\right.\\
  &=& \left\{\begin{array}{lcl}
           \hat{w}\sigma_i U^i T &\mb{ for } &  w_2\al_i\in \nrW\setminus\We_{\Th\cup\Th^\bot}(\Th\cup\Th^\bot)\cup\We_\Th\Th\\
           U_{w_1 w_2\al_i}\hat{w}\sigma_i U^i T &\mb{ for }& w_2\al_i\in\We_{\Th^\bot}\Th^\bot
      \end{array}\right.\;\;.
\end{eqnarray*}
Since $w_1\in\We^{\Th\cup\Th^\bot}$ we have $w_1 w_2\al_i\in (\rW)^\pm$ if $w_2\al_i\in\We_{\Th^\bot}\Th^\bot \cap (\rW)^\pm$. We conclude:
\begin{eqnarray*}
   B \hat{w} B \sigma_i B    &=&  B\hat{w}\sigma_i B  \quad\mb{ for }\quad w_2\al_i\in 
                                 \left(\nrW\setminus \We_{\Th\cup\Th^\bot}(\Th\cup\Th^\bot) \right) \cup 
                                  (\We_{\Th^\bot}\Th^\bot \cap \prW) \cup\We_\Th\Th\;\;.\\
   B^- \hat{w} B \sigma_i B  &=&  B^-\hat{w}\sigma_i B    \quad\mb{ for }\quad w_2\al_i\in\left(\nrW\setminus \We_\Th \Th\right)\cup\We_\Th\Th\;\;.
\end{eqnarray*}
{\bf b)} We have
\begin{eqnarray*}
  \hat{w} B \sigma_i &=& w_1\ve{R(\Th)} w_2 U_i U^i T \sigma_i   \;\,=\;\, w_1\ve{R(\Th)} w_2 U_i \sigma_i U^i T \\
                     &=&  w_1\ve{R(\Th)} w_2 \sigma_i(\sigma_i^{-1} U_i \sigma_i) U^i T
                     \;\,\subseteq\;\, 
                            w_1\ve{R(\Th)} w_2 \sigma_i (U_i\cup U_i\sigma_i U_i) U^i T \\
                     &=&  \hat{w}\sigma_i B \;\,\cup\;\, w_1\ve{R(\Th)} U_{w_2 \sigma_i\al_i} w_2  B \\
                     &=&  \hat{w}\sigma_i B \;\,\cup\;\, \left\{ \begin{array}{lcl}
                          w_1\ve{R(\Th)}w_2 B & for &  w_2\al_i \in \prW\setminus \We_{\Th\cup\Th^\bot}(\Th\cup\Th^\bot)\\
                          w_1 U_{-w_2\al_i}\ve{R(\Th)} w_2 B& for &  w_2\al_i \in \We_{\Th^\bot}\Th^\bot
                          \end{array}\right.\\
                     &=&  \hat{w}\sigma_i B \;\,\cup\;\, \left\{ \begin{array}{lcl}
                          \hat{w} B & for &  w_2\al_i \in \prW\setminus \We_{\Th\cup\Th^\bot}(\Th\cup\Th^\bot)\\
                          U_{-w_1 w_2\al_i}\hat{w}B & for &  w_2\al_i \in \We_{\Th^\bot}\Th^\bot
                          \end{array}\right.\;\;.\\                                             
\end{eqnarray*}
Here $-w_1w_2\al_i\in (\rW)^\mp$ if $w_2\al_i\in\We_{\Th^\bot}\Th^\bot \cap (\rW)^\pm$. Note also that 
$\sigma_i^{-1} U_i\sigma_i$ contains elements of $U_i$ and $U_i\sigma_i U_i$. We conclude:
\begin{eqnarray*}
   B \hat{w} B \sigma_i B      &=& B\hat{w}\sigma_i B \cup  B\hat{w}B \quad\mb{ for }\quad w_2\al_i\in 
                                 \left(\prW\setminus \We_{\Th\cup\Th^\bot}(\Th\cup\Th^\bot)\right) \cup 
                                   (\We_{\Th^\bot}\Th^\bot \cap \nrW) \;\;.\\
   B^- \hat{w} B \sigma_i B    &=&  B^-\hat{w}\sigma_i B \cup  B^-\hat{w}B    \quad\mb{ for }\quad w_2\al_i\in\prW\setminus \We_\Th \Th\;\;.
\end{eqnarray*}
Due to Part 2) of the last theorem, from a) and b) follows:
\begin{eqnarray*}
   B \hat{w} B \sigma_i B &=& \left\{\begin{array}{lcl}
                                           B\hat{w}B                             &\mb{if}&   l_{++}(\hat{w}\sigma_i)=l_{++}(\hat{w}) \\
                                           B\hat{w}\sigma_i B                    &\mb{if}&   l_{++}(\hat{w}\sigma_i)=l_{++}(\hat{w})+1 \\
                                           B\hat{w}\sigma_i B \;\dot{\cup}\; B\hat{w}B     &\mb{if}&   l_{++}(\hat{w}\sigma_i)=l_{++}(\hat{w})-1 
                                          \end{array}\right.\;\;.\\
   B^- \hat{w} B \sigma_i B &=& \left\{\begin{array}{lcl}
                                           B^-\hat{w}B                               &\mb{if}&   l_{--}(\hat{w}\sigma_i)=l_{--}(\hat{w}) \\
                                           B^-\hat{w}\sigma_i B                      &\mb{if}&   l_{--}(\hat{w}\sigma_i)=l_{--}(\hat{w})-1 \\
                                           B^-\hat{w}\sigma_i B \;\dot{\cup}\; B^-\hat{w}B     &\mb{if}&   l_{--}(\hat{w}\sigma_i)=l_{--}(\hat{w})+1 
                                          \end{array}\right.\;\;.
\end{eqnarray*}
The remaining cases of the theorem follow by applying the Chevalley involution $*:\GD\to\GD$, together with 
Remark (3) following Definition \ref{BO1} and Remark (1) following Definition \ref{L1}.\\
\End
%
%
%
%
%
%
%
%
%
{\bf Acknowledgment.} I would like to thank the Deutsche Forschungsgemeinschaft for providing my main financial 
support, when most of this article has been elaborated. I also would like to thank the TMR-Program 
ERB FMRX-CT97-0100 ``Algebraic Lie Theory'' for providing some financial support for traveling.
%
%
%
%
%

%
%
%
%
\end{document}